# The Derivation of Markov Chain Properties using Generalized Matrix Inverses

Jeffrey J. Hunter
School of Computing and Mathematical Sciences
Auckland University of Technology
Auckland, New Zealand

**Abstract**
The analysis of many problems of interest associated with Markov chains, e.g. stationary distributions, moments of first passage time distributions and moments of occupation time random variables, involves the solution of a system of linear equations involving *I* - *P*, where *P* is the transition matrix of a finite, irreducible, discrete time Markov chain. Generalized inverses play an important role in the solution of such singular sets of equations. In this presentation we survey the application of generalized matrix inverses to the aforementioned problems.

---

## 1. Introduction

In many instances of stochastic modelling a Markov chain is typically present, either on its own or embedded within a more general process, and one wishes to determine its key properties such as the stationary distribution, the mean first passage time between states (and higher order moments) and expected occupation times in particular states given the Markov chain starts in a specified state.

We shall see that these key properties invariably involve the solution of sets of constrained linear equations. The solution to these equations can be effected by using generalized inverses.

We first introduce some basic notation for Markov chains, explore the equations to be solved, look at the general solution of systems of linear equations, using generalized matrix inverses. Following a discussion of generalized inverses of Markovian kernels, we provide a systematic presentation of solving for the key properties.

## 2. Properties of Markov chains

Let $\{X_n\}$ be a Markov chain with a finite state space $S = \{1,2,...,m\}$ and transition matrix $P = [p_{ij}]$ where $p_{ij} = P[X_n = j \mid X_{n-1} = i]$ for all $i, j \in S$, $n \geq 1$.

### 2.1 Stationary distributions

It is well known that if the Markov chain is irreducible (every state can be reached from every state), then a stationary distribution $\{\pi_i, i \in S\}$ exists, which is the limiting distribution in the case of aperiodic Markov chains.

The stationary probability vector $\boldsymbol{\pi}^T = (\pi_1, \pi_2, ..., \pi_m)$ for a finite irreducible Markov chain satisfies the set of linear equations

$$\pi_j = \sum_{i=1}^{m} \pi_i p_{ij}, \text{ with the constraint that } \sum_{i=1}^{m} \pi_i = 1. \tag{2.1}$$

or, in matrix form,

$$\boldsymbol{\pi}^T (I - P) = \boldsymbol{0}^T \text{ with } \boldsymbol{\pi}^T \boldsymbol{e} = 1, \tag{2.2}$$

where $\boldsymbol{e}^T = (1,1,...,1)$.

Equation (2.2) is an equation of type $XB = C$ with $X = \boldsymbol{\pi}^T$, $B = I - P$, $C = \boldsymbol{0}^T$.

**2.2 Moments of first passage times**

Let $T_{ij}$ be the number of trials for a Markov chain $\{X_n\}$ to make a first passage from state $i$ to state $j$ (first return if $i = j$) i.e. $T_{ij} = \min\{n\ (\geq 1) : X_n = j \text{ given } X_0 = i\}$. Under the assumption of irreducibility for finite Markov chains, the $T_{ij}$ are proper random variables and hence the moments $m_{ij}^{(k)} = E[T_{ij}^k \mid X_0 = i]$ $(k \geq 1)$ are well defined and finite for all $i, j \in S$, (Neuts, 1973), (Hunter, 1983)). It is well known that, for all $i, j \in S$, the first moments of the first passage times $m_{ij} = m_{ij}^{(1)}$ satisfy the set of linear equations ( Hunter, 1983)

$$m_{ij} = p_{ij} + \sum_{k \neq j} p_{ik}(m_{kj} + 1)$$

i.e.

$$m_{ij} = 1 + \sum_{k \neq j} p_{ik} m_{kj}. \tag{2.3}$$

Define $M = [m_{ij}]$, $M_d = [\delta_{ij} m_{ij}]$, (where $\delta_{ij} = 1$ if $i = j$ and 0 otherwise), i.e. $M_d$ is a diagonal matrix with elements the diagonal elements of $M$ and let $E = \boldsymbol{e}\boldsymbol{e}^T = [1]$, (i.e. all the elements of E are unity). Equations (2.3) can be expressed as

$$(I - P)M = E - PM_d, \tag{2.4}$$

It is also well known (Feller, 1950), (Kemeny and Snell, 1960), (Hunter, 1983) that

$$m_{ii} = \frac{1}{\pi_i}.$$

This implies

$$M_d = D = \left[\delta_{ij} / \pi_i\right] = (\Pi_d)^{-1}, \tag{2.5}$$

where $\Pi = \boldsymbol{e}\boldsymbol{\pi}^T$, a matrix with all rows identical to the stationary probability vector.

(When the Markov chain is aperiodic, $\Pi = \lim_{n \to \infty} [p_{ij}^{(n)}] = [\pi_j]$, i.e. the rows of $\Pi$ are the limiting distribution of the Markov chain.)

To solve for $M$ we have from equations (2.4) and (2.5), a matrix equation of the form $AX = C$ where $A = I - P$, $X = M$ and $C = E - PD$, which is fully specified.

Further it can be shown (Kemeny and Snell, 1960), (Hunter, 1983) that $M^{(2)}$ satisfies the matrix equation

$$(I - P)M^{(2)} = E + 2P(M - M_d) - PM_d^{(2)}. \tag{2.6}$$

We show in Section 7 that $M_d^{(2)}$ can be expressed in terms of the matrices, $M$, $\Pi$ and $D$,

$$M_d^{(2)} = 2D(\Pi M)_d - D.$$

This also leads to a matrix equation for $M^{(2)}$ of the form $AX = C$ where once again $A = I - P$ with $X = M^{(2)}$ and $C$ a predetermined matrix.

Matrix equations for the higher moments $M^{(k)}$ can also be developed, (see Hunter, 1990).

**2.3 Occupation times**

Another application arises in examining the asymptotic behaviour of the number of times particular states are entered. Given a Markov chain $\{X_n\}$, let

$M_{ij}^{(n)} = \{$Number of $k$ $(0 \le k \le n)$ such that $X_k = j$ given $X_0 = i\}$, then

$$\left[EM_{ij}^{(n)}\right] = \left[\sum_{k=0}^{n} p_{ij}^{(k)}\right] = \sum_{k=0}^{n} P^k .$$

Let $A_n = \sum_{k=0}^{n-1} P^k$ then

$$(I - P)A_n = I - P^n \text{ or } A_n(I - P) = I - P^n, \tag{2.7}$$

subject to $A_n\Pi = \Pi A_n = n\Pi$.

Equations (2.7) with unknown $X = A_n$ are of the form $AX = C$ or $XB = C$, where $A = B = I - P$.

**2.4 Further extensions**

In this workshop we consider only the properties of Markov chains. However many of the results that we develop can be extended to examine the key related properties of Markov renewal processes, semi-Markov processes and Markov chains in continuous time (where a Markov chain is embedded within the processes of interest.) For further information consult (Hunter, 1982).

## 3. Generalised matrix inverses

A generalized matrix inverse (g-inverse) of a matrix $A$ is any matrix $A^-$ such that $AA^-A = A$.
If $A$ is non-singluar then $A^-$ is $A^{-1}$, the inverse of $A$, and is unique. If $A$ is singular, $A^-$ is not unique.
Matrices $A^-$ are called "one condition" g-inverses or "equation solving" g-inverses because, as we shall see, of their use in solving systems of linear equations.

By imposing additional conditions we can obtain various types of multi-condition g-inverses. Consider real conformable matrices X such that:

(Condition 1) $AXA = A$.
(Condition 2) $XAX = X$.
(Condition 3) $(AX)^T = AX$.
(Condition 4) $(XA)^T = XA$.
(Condition 5) For square matrices $AX = XA$.

Let $A^{(i,j,\ldots,l)}$ be any matrix that satisfies conditions (*i*), (*j*), ..., (*l*) of the above itemised conditions. $A^{(i,j,\ldots,l)}$ is called an - (*i*, *j*, ..., *l*) g-inverse of $A$, under the assumption that condition (1) is always included.

Particular special cases include $A^{(1,2)}$ - a "pseudo-inverse" (Rao, 1955), or "reciprocal inverse" (Bjerhammar, 1951), or a "reflexive" g-inverse (Rhode, 1964); $A^{(1,3)}$ - a "least squares" g-inverse; $A^{(1,4)}$ - a "minimum norm" g-inverse; $A^{(1,2,4)}$ - a "weak generalized inverse" (Goldman & Zelen, 1964); $A^{(1,2,3,4)}$ - the "Moore-Penrose" g-inverse, (Moore, 1920), (Penrose, 1955); and $A^{(1,2,5)}$ - the "group inverse" (Erdelyi, 1967) which exists and is unique if $\text{rank}(A) = \text{rank}(A^2)$. Note that, except for $A^{(1,2,3,4)}$ and $A^{(1,2,5)}$, g-inverses are in general not unique.

## 4. Solving Systems of Linear Equations

Generalized inverses play a major role in solving systems of linear equations.

**Theorem 4.1**: (Penrose, 1955), (Rao, 1955) *A necessary and sufficient condition for $AXB = C$ to have a solution is that*

$$AA^{-}CB^{-}B = C. \tag{4.1}$$

*If this consistency condition is satisfied the general solution is given by*

$$X = A^{-}CB^{-} + W - A^{-}AWBB^{-}, \tag{4.2}$$

*where W is an arbitrary matrix*, or

$$X = A^{-}CB^{-} + (I - A^{-}A)U + V(I - BB^{-}), \tag{4.3}$$

*where U and V are arbitrary matrices.*

**Proof**: If $AXB = C$ is consistent there exists an $X_0$ such that $C = AX_0B$. Thus $AA^{-}CB^{-}B = AA^{-}AX_0BB^{-}B = AX_0B = C.$

If $C = AA^{-}CB^{-}B$ then $X = A^{-}CB^{-}$ is a particular solution.

To establish the general solution (equation (4.2), first note that for every $W$, $A^{-}CB^{-} + W - A^{-}AWBB^{-}$ is a solution of $AXB = C$. Further, given a particular solution $X_0$ there exists a $W_0$ such that $X_0 = A^{-}CB^{-} + W_0 - A^{-}AW_0BB^{-}$. In particular, if we take $W_0 = X_0$ then $AX_0B = AA^{-}CB^{-}B + AX_0B - AA^{-}AX_0BB^{-}B = C + AX_0B - AX_0B = C$.

The equivalence of the alternative general solution, equation (4.3), follows by taking $W = (I - A^{-}A)U + V(I - BB^{-})$ or, conversely, by taking $U = V = \frac{1}{2}(W + A^{-}AW + WBB^{-})$.

Special cases of interest are the following:

(i) $XB = C$ has the general solution $X = CB^{-} + W(I - BB^{-})$ where $W$ is an arbitrary matrix, provided $CB^{-}B = C$. (4.4)

(ii) $AX = C$ has the general solution $X = A^{-}C + (I - A^{-}A)W$, where $W$ is an arbitrary matrix, provided $AA^{-}C = C$. (4.5)

(iii) $AXA = A$ has the general solution $X = A^{-}AA^{-} + W - A^{-}AWAA^{-}$, where $W$ is an arbitrary matrix. (4.6)

**Theorem 4.2**: *If $A^{-}$ is any g-inverse of A then all g-inverses of A can be characterised as members of the following equivalent sets:*

$$A\{1\} = \{A^{-}AA^{-} + W - A^{-}AWAA^{-},\ W \text{ arbitrary}\}. \tag{4.7}$$

$$A\{1\} = \{A^{-}AA^{-} + (I - A^{-}A)U + V(I - AA^{-}),\ U,\ V \text{ arbitrary}\}. \tag{4.8}$$

$$A\{1\} = \{A^{-} + H - A^{-}AHAA^{-},\ H \text{ arbitrary}\}. \tag{4.9}$$

$$A\{1\} = \{A^- + (I - A^- A)F + G(I - AA^-),\ F,\ G \text{ arbitrary}\}. \quad (4.10)$$

**Proof**: The first two characterisations (4.7) and (4.8) follow from (4.2) and (4.3) with $C = B = A$. Furthermore by taking $W = H + A^-$, (4.7) and (4.9) are equivalent as are (4.8) and (4.10) with $U \equiv F + \frac{1}{2}A^-$ and $V \equiv G + \frac{1}{2}A^-$.

## 5. Generalized inverses of Markovian kernels

**Theorem 5.1**: (Hunter, 1982) *Let $P$ be the transition matrix of a finite irreducible Markov chain with m states and stationary probability vector $\boldsymbol{\pi}^T = (\pi_1, \pi_2, \ldots, \pi_m)$.*
*Let $\boldsymbol{e}^T = (1, 1, ..., 1)$ and $\boldsymbol{t}$ and $\boldsymbol{u}$ be any vectors.*
(a) *$I - P + \boldsymbol{tu}^T$ is non-singular if and only if $\boldsymbol{\pi}^T\boldsymbol{t} \neq 0$ and $\boldsymbol{u}^T\boldsymbol{e} \neq 0$.*
(b) *If $\boldsymbol{\pi}^T\boldsymbol{t} \neq 0$ and $\boldsymbol{u}^T\boldsymbol{e} \neq 0$ then $[I - P + \boldsymbol{tu}^T]^{-1}$ is a g-inverse of $I - P$.*

**Proof:** (a) For any matrix $X$, $\det(X + \boldsymbol{tu}^T) = \det(X) + \boldsymbol{u}^T[\text{adj } X]\boldsymbol{t}$.
Thus, taking $X = I - P$, a singular matrix, $\det(I - P + \boldsymbol{tu}^T) = \boldsymbol{u}^T[\text{adj } (I - P)]\boldsymbol{t}$.
If $A = \text{adj } (I - P)]$ then $A(I - P) = (I - P)A = \det(I - P)I = 0$, so that

$$A = AP = PA. \quad (5.1)$$

If $P$ is irreducible, then any matrix satisfying equation (5.1) is a multiple of $\Pi = \boldsymbol{e\pi}^T$ implying that $A = \text{adj } (I - P)] = k\boldsymbol{e\pi}^T$. Since $P$ is irreducible, its eigenvalues $\lambda_1, \lambda_2, ..., \lambda_m$ are such that $\lambda_1 = 1$ is the only eigenvalue equal to 1. Consequently $\text{tr}(\text{adj } (I - P)) = \prod_{j=2}^{m}(1 - \lambda_j) \neq 0$. But $\text{tr}(\text{adj } (I - P)) = ktr(\boldsymbol{e\pi}^T) = k(\sum_{i=1}^{m}\pi_i) = k \neq 0$. Thus $\det(I - P + \boldsymbol{tu}^T) = \boldsymbol{u}^T A\boldsymbol{t} = k\boldsymbol{u}^T\boldsymbol{e\pi}^T\boldsymbol{t} = k(\boldsymbol{u}^T\boldsymbol{e})(\boldsymbol{\pi}^T\boldsymbol{t}) \neq 0$ establishing the non-singularity.

(b) Observe that $(I - P + \boldsymbol{tu}^T)(I - P + \boldsymbol{tu}^T)^{-1} = I$ so that

$$(I - P)(I - P + \boldsymbol{tu}^T)^{-1} = I - \boldsymbol{tu}^T(I - P + \boldsymbol{tu}^T)^{-1}.$$

Now $\quad \boldsymbol{\pi}^T(I - P + \boldsymbol{tu}^T) = \boldsymbol{\pi}^T(I - P) + \boldsymbol{\pi}^T\boldsymbol{tu}^T = (\boldsymbol{\pi}^T\boldsymbol{t})\boldsymbol{u}^T$,

so that $\quad \boldsymbol{\pi}^T = (\boldsymbol{\pi}^T\boldsymbol{t})\boldsymbol{u}^T(I - P + \boldsymbol{tu}^T)^{-1} \Rightarrow \dfrac{\boldsymbol{\pi}^T}{\boldsymbol{\pi}^T\boldsymbol{t}} = \boldsymbol{u}^T(I - P + \boldsymbol{tu}^T)^{-1}$.

Further, $(I - P)(I - P + \boldsymbol{tu}^T)^{-1} = I - \boldsymbol{tu}^T(I - P + \boldsymbol{tu}^T)^{-1} = I - \dfrac{\boldsymbol{t\pi}^T}{\boldsymbol{\pi}^T\boldsymbol{t}}$,

implying $(I - P)(I - P + \boldsymbol{tu}^T)^{-1}(I - P) = I - P - \dfrac{\boldsymbol{t\pi}^T(I - P)}{\boldsymbol{\pi}^T\boldsymbol{t}} = I - P$

and hence establishing that $(I - P + \boldsymbol{tu}^T)^{-1}$ is a "one condition" g-inverse of $I - P$.
[Alternatively note that $(I - P + \boldsymbol{tu}^T)^{-1}(I - P + \boldsymbol{tu}^T) = I$ so that

$$(I - P + \boldsymbol{tu}^T)^{-1}(I - P) = I - (I - P + \boldsymbol{tu}^T)^{-1}\boldsymbol{tu}^T.$$

Now $(I - P + \boldsymbol{tu}^T)\boldsymbol{e} = (I - P)\boldsymbol{e} + \boldsymbol{tu}^T\boldsymbol{e} = \boldsymbol{t}(\boldsymbol{u}^T\boldsymbol{e})$, so that $\boldsymbol{e} = (I - P + \boldsymbol{tu}^T)^{-1}\boldsymbol{t}(\boldsymbol{u}^T\boldsymbol{e})$

$\Rightarrow \frac{e}{u^T e} = (I - P + tu^T)^{-1} t$ . Further $(I - P + tu^T)^{-1}(I - P) = I - \frac{eu^T}{u^T e}$,

so that $(I - P)(I - P + tu^T)^{-1}(I - P) = (I - P)(I - \frac{eu^T}{u^T e}) = I - P - (I - P)\frac{eu^T}{u^T e} = I - P$,

providing an alternative proof of (b).]

We have in the process of the above proof established the following useful results:

(a) $$u^T [I - P + tu^T]^{-1} = \frac{\pi^T}{\pi^T t}. \qquad (5.2)$$

(b) $$[I - P + tu^T]^{-1} t = \frac{e}{u^T e}. \qquad (5.3)$$

**Theorem 5.2**: *Under the conditions of Theorem 5.1, all "one condition" g-inverses of I – P can be expressed by any of the following equivalent forms:*

(i) $A^{(1)} = [I - P + tu^T]^{-1} + \frac{eu^T H}{u^T e} + \frac{Ht\pi^T}{\pi^T t} - \frac{eu^T Ht\pi^T}{(u^T e)(\pi^T t)}$ *for arbitrary matrix H.* (5.4)

(ii) $A^{(1)} = [I - P + tu^T]^{-1} + \frac{eu^T F}{u^T e} + \frac{Gt\pi^T}{(\pi^T t)}$ *for arbitrary matrices F and G.* (5.5)

(iiii) $A^{(1)} = [I - P + tu^T]^{-1} + ef^T + g\pi^T$ *for arbitrary vectors **f** and **g**.* (5.6)

**Proof**: Use Theorem 4.2, and subsidiary results used in establishing Theorem 5.1 including equations (5.2) and (5.3).

The presenter published a paper (Hunter, 1988) on various characterisations of generalized inverses associated with Markovian kernels. A systematic investigation into various multi-condition generalized inverses based upon the results of Theorem 5.1 including a fully efficient characterization of a "one condtion" g-inverse of *I – P* was given. In particular it was shown that:

**Theorem 5.3**: *Under the conditions of Theorem 5.1, all "one condition" g-inverses of I – P can be expressed as*

$$A^{(1)} = [I - P + tu^T]^{-1} + ey_m^T + Bz_{m-1}\pi^T, \qquad (5.7)$$

*where if* $u^T = (u_{m-1}^T, u_m)$, $B = \begin{bmatrix} I - \left(\frac{1}{u^T e}\right) eu_{m-1}^T \\ \cdots\cdots\cdots \\ -\left(\frac{1}{u^T e}\right) eu_{m-1}^T \end{bmatrix}$, *an m × (m – 1) matrix,* $y_{m-1}^T$ *and* $z_{m-1}^T$ *are arbitrary vectors with respectively m and m – 1 elements.*

The proof (omitted) was based upon results of (Ben-Israel & Greville , 1974).

Observe that $P$ is $m \times m$ and $\boldsymbol{t}$ and $\boldsymbol{u}$ are given, then the expression in equation (5.6) has $2m$ arbitrary elements. The characterization given by equation (5.7), with specified $\boldsymbol{t}$ and $\boldsymbol{u}$, however requires only $2m - 1$ arbitrary elements. This representation is fully efficient in that given any generalized inverse, of $I - P$ of the form given by equation (5.7) with prescribed $\boldsymbol{t}$ and $\boldsymbol{u}$,

$$\boldsymbol{y}_m^T = \frac{\boldsymbol{u}^T}{\boldsymbol{u}^T\boldsymbol{e}}\left\{A^{(1)} - (I - P + \boldsymbol{t}\boldsymbol{u}^T)^{-1}\right\},\ \boldsymbol{z}_{m-1}^T = \left[\ I \ \vdots \ -\boldsymbol{e}\ \right]\left\{A^{(1)} - (I - P + \boldsymbol{t}\boldsymbol{u}^T)^{-1}\right\}\frac{\boldsymbol{t}}{\boldsymbol{\pi}^T\boldsymbol{t}}.$$

There are some disadvantages in that representation (5.7) requires preselected $\boldsymbol{t}$, $\boldsymbol{u}$.

However (see Hunter, 1988) we can convert from a form of one generalized inverse of $I - P$ with given $\boldsymbol{t}_1$, $\boldsymbol{u}_1$ to another equivalent form with $\boldsymbol{t}_2$, $\boldsymbol{u}_2$ (with $\boldsymbol{\pi}^T\boldsymbol{t}_i \neq 0$ and $\boldsymbol{u}_i^T\boldsymbol{e} \neq 0$, $(i = 1, 2)$):

$$\left[I - P + \boldsymbol{t}_2\boldsymbol{u}_2^T\right]^{-1} = \left[I - \frac{\boldsymbol{e}\boldsymbol{u}_2^T}{\boldsymbol{u}_2^T\boldsymbol{e}}\right]\left[I - P + \boldsymbol{t}_1\boldsymbol{u}_1^T\right]^{-1}\left[I - \frac{\boldsymbol{t}_2\boldsymbol{\pi}^T}{\boldsymbol{\pi}^T\boldsymbol{t}_2}\right] + \frac{\boldsymbol{e}\boldsymbol{\pi}^T}{(\boldsymbol{\pi}^T\boldsymbol{t}_2)(\boldsymbol{u}_2^T\boldsymbol{e})}.$$

The proof (omitted) uses the result that $(X + \boldsymbol{a}\boldsymbol{b}^T)^{-1} = X^{-1} - \dfrac{(X^{-1}\boldsymbol{a})(\boldsymbol{b}^T X^{-1})}{1 + \boldsymbol{b}^T X^{-1}\boldsymbol{a}}$.

Further, (Hunter, 1990), given any g-inverse $G$ of $I - P$ and $\boldsymbol{t}$, $\boldsymbol{u}$, with $\boldsymbol{\pi}^T\boldsymbol{t} \neq 0$ and $\boldsymbol{u}^T\boldsymbol{e} \neq 0$, we can compute $\left[I - P + \boldsymbol{t}\boldsymbol{u}^T\right]^{-1}$:

$$\left[I - P + \boldsymbol{t}\boldsymbol{u}^T\right]^{-1} = \left[I - \frac{\boldsymbol{e}\boldsymbol{u}^T}{\boldsymbol{u}^T\boldsymbol{e}}\right]G\left[I - \frac{\boldsymbol{t}\boldsymbol{\pi}^T}{\boldsymbol{\pi}^T\boldsymbol{t}}\right] + \frac{\boldsymbol{e}\boldsymbol{\pi}^T}{(\boldsymbol{\pi}^T\boldsymbol{t})(\boldsymbol{u}^T\boldsymbol{e})}.$$

A further subsidiary result, (Hunter, 1988), is that

$$A^{(1)} = [I - P + \delta\boldsymbol{t}\boldsymbol{u}^T]^{-1} - \frac{\boldsymbol{e}\boldsymbol{\pi}^T}{\delta(\boldsymbol{\pi}^T\boldsymbol{t})(\boldsymbol{u}^T\boldsymbol{e})} \text{ does not depend on } \delta \neq 0. \tag{5.8}$$

In (Hunter, 1988) conditions for the various multi-condition g-inverses of $I - P$ were derived:

**Theorem 5.4**: *Under the assumption that* $\boldsymbol{\pi}^T\boldsymbol{t} \neq 0$ *and* $\boldsymbol{u}^T\boldsymbol{e} \neq 0$, *the general one condition g-inverse of the form* $G = [I - P + \boldsymbol{t}\boldsymbol{u}^T]^{-1} + \boldsymbol{e}\boldsymbol{f}^T + \boldsymbol{g}\boldsymbol{\pi}^T$ *is a member of the following multi-condition g-inverse families*:

(1) $G \in A\{1,2\} \Leftrightarrow \dfrac{1}{(\boldsymbol{u}^T\boldsymbol{e})(\boldsymbol{\pi}^T\boldsymbol{t})} + \dfrac{\boldsymbol{f}^T\boldsymbol{t}}{\boldsymbol{\pi}^T\boldsymbol{t}} + \dfrac{\boldsymbol{u}^T\boldsymbol{g}}{\boldsymbol{u}^T\boldsymbol{e}} = \boldsymbol{f}^T(I - P)\boldsymbol{g}.$

(2) $G \in A\{1,3\} \Leftrightarrow \dfrac{\boldsymbol{t}}{\boldsymbol{\pi}^T\boldsymbol{t}} - (I - P)\boldsymbol{g} = \dfrac{\boldsymbol{\pi}}{\boldsymbol{\pi}^T\boldsymbol{\pi}}.$

(3) $G \in A\{1,4\} \Leftrightarrow \dfrac{\boldsymbol{u}^T}{\boldsymbol{u}^T\boldsymbol{e}} - \boldsymbol{f}^T(I - P) = \dfrac{\boldsymbol{e}^T}{\boldsymbol{e}^T\boldsymbol{e}}.$

(4) $G \in A\{1,5\} \Leftrightarrow \dfrac{\boldsymbol{t}}{\boldsymbol{\pi}^T\boldsymbol{t}} - (I - P)\boldsymbol{g} = \boldsymbol{e}$ and $\dfrac{\boldsymbol{u}^T}{\boldsymbol{u}^T\boldsymbol{e}} - \boldsymbol{f}^T(I - P) = \boldsymbol{\pi}^T.$

The big disadvantage of these above results and conditions is that they are based on the representation for a typical one condition g-inverse $A^{(1)}$ of the form $A^{(1)} = [I - P + \boldsymbol{tu}^T]^{-1} + \boldsymbol{ef}^T + \boldsymbol{g\pi}^T$, as given by equation (5.6) and consequently requires knowledge of the $\boldsymbol{t}$, $\boldsymbol{u}$, $\boldsymbol{f}$ and $\boldsymbol{g}$.

In a subsequent paper (Hunter, 1990) parametric forms of g-inverses of $I - P$ were obtained. As expected, these forms still require $2m - 1$ independent parameters but they have the feature that they uniquely determine the particular g-inverse and its g-inverse. The key result is the following.

**Theorem 5.5**: *Given any g-inverse, G, of I – P there exists unique parameters* $\boldsymbol{\alpha}$, $\boldsymbol{\beta}$, $\gamma$ *with the property that*

$$G = [I - P + \boldsymbol{\alpha\beta}^T]^{-1} + \gamma \boldsymbol{e\pi}^T, \tag{5.9}$$

where $\boldsymbol{\pi}^T\boldsymbol{\alpha} = 1$, $\boldsymbol{\beta}^T\boldsymbol{e} = 1$ and $\gamma + 1 = \boldsymbol{\pi}^T G\boldsymbol{\alpha} = \boldsymbol{\beta}^T G\boldsymbol{e} = \boldsymbol{\beta}^T G\boldsymbol{\alpha}$.

**Proof**: Let $G$ be any g-inverse of $I - P$, so that in general form, from equation (5.6), $G$ can be expressed as

$$G = [I - P + \boldsymbol{tu}^T]^{-1} + \boldsymbol{ef}^T + \boldsymbol{g\pi}^T.$$

Now define the following matrices and vectors

$$A \equiv I - (I - P)G \tag{5.10}$$

$$B \equiv I - G(I - P), \tag{5.11}$$

$$\boldsymbol{\alpha} \equiv \frac{\boldsymbol{t}}{\boldsymbol{\pi}^T\boldsymbol{t}} - (I - P)\boldsymbol{g}, \tag{5.12}$$

$$\boldsymbol{\beta}^T \equiv \frac{\boldsymbol{u}^T}{\boldsymbol{u}^T\boldsymbol{e}} - \boldsymbol{f}^T(I - P). \tag{5.13}$$

First note that $\boldsymbol{\pi}^T\boldsymbol{\alpha} = 1$ and $\boldsymbol{\beta}^T\boldsymbol{e} = 1$. (5.14)

Further, from the proof of Theorem 5.1,

$$(I - P)(I - P + \boldsymbol{tu}^T)^{-1} = I - \frac{\boldsymbol{t\pi}^T}{\boldsymbol{\pi}^T\boldsymbol{t}},$$

and

$$(I - P + \boldsymbol{tu}^T)^{-1}(I - P) = I - \frac{\boldsymbol{eu}^T}{\boldsymbol{u}^T\boldsymbol{e}},$$

leading to the observations that

$$A = \boldsymbol{\alpha\pi}^T \text{ and } B = \boldsymbol{e\beta}^T. \tag{5.15}$$

Further $\boldsymbol{\alpha}$ and $\boldsymbol{\beta}$ can be recovered from $A$ and $B$ simply as

$$\boldsymbol{\alpha} = A\boldsymbol{e} \ (\neq 0) \text{ and } \boldsymbol{\beta}^T = \boldsymbol{\pi}^T B \ (\neq 0). \tag{5.16}$$

Note also that $A\boldsymbol{\alpha} = \boldsymbol{\alpha}$, $\boldsymbol{\pi}^T A = \boldsymbol{\pi}^T$, $\boldsymbol{\beta}^T B = \boldsymbol{\beta}^T$ and $B\boldsymbol{e} = \boldsymbol{e}$.

As noted earlier there is considerable flexibility in how we choose $\boldsymbol{t}$ and $\boldsymbol{u}$. It would be advantageous if, for a general g-inverse of the form given by equation (5.6), we could choose $\boldsymbol{t}$ and $\boldsymbol{u}$ so that we can easily determine explicit expressions for the vectors $\boldsymbol{f}$ and $\boldsymbol{g}$. Let us take $\boldsymbol{t} = A\boldsymbol{e} = \boldsymbol{\alpha}$ and $\boldsymbol{u}^T = \boldsymbol{\pi}^T B = \boldsymbol{\beta}^T$.

Then, from the definition (5.12),

$$\boldsymbol{\alpha} \equiv \frac{\boldsymbol{\alpha}}{\boldsymbol{\pi}^T\boldsymbol{\alpha}} - (I-P)\boldsymbol{g} = \boldsymbol{\alpha} - (I-P)\boldsymbol{g} \Rightarrow (I-P)\boldsymbol{g} = \boldsymbol{0} \Rightarrow \boldsymbol{g} = g\boldsymbol{e} \text{ for some } g,$$

and, from definition (5.13),

$$\boldsymbol{\beta}^T = \frac{\boldsymbol{\beta}^T}{\boldsymbol{\beta}^T\boldsymbol{e}} - \boldsymbol{f}^T(I-P) = \boldsymbol{\beta}^T - \boldsymbol{f}^T(I-P) \Rightarrow \boldsymbol{f}^T(I-P) = \boldsymbol{0}^T \Rightarrow \boldsymbol{f}^T = f\boldsymbol{\pi}^T \text{ for some } f.$$

Thus with this particular choice of $\boldsymbol{t}$ and $\boldsymbol{u}$, $G$ takes the form given by equation (5.9) with $\gamma = f + g$. To determine expressions for $\gamma$ note that from equations (5.2) and (5.3) (respectively),

$$\boldsymbol{\beta}^T[I-P+\boldsymbol{\alpha}\boldsymbol{\beta}^T]^{-1} = \frac{\boldsymbol{\pi}^T}{\boldsymbol{\pi}^T\boldsymbol{\alpha}} = \boldsymbol{\pi}^T \text{ and } [I-P+\boldsymbol{\alpha}\boldsymbol{\beta}^T]^{-1}\boldsymbol{\alpha} = \frac{\boldsymbol{e}}{\boldsymbol{\beta}^T\boldsymbol{e}} = \boldsymbol{e}.$$

Thus, with $G$ given by equation (5.9),

$$\boldsymbol{\beta}^T G = \boldsymbol{\beta}^T[I-P+\boldsymbol{\alpha}\boldsymbol{\beta}^T]^{-1} + \gamma\boldsymbol{\beta}^T\boldsymbol{e}\boldsymbol{\pi}^T = (\gamma+1)\boldsymbol{\pi}^T \Rightarrow \boldsymbol{\beta}^T G\boldsymbol{e} = \boldsymbol{\beta}^T G\boldsymbol{\alpha} = \gamma + 1,$$

and $$G\boldsymbol{\alpha} = [I-P+\boldsymbol{\alpha}\boldsymbol{\beta}^T]^{-1}\boldsymbol{\alpha} + \gamma\boldsymbol{e}\boldsymbol{\pi}^T\boldsymbol{\alpha} = (\gamma+1)\boldsymbol{e} \Rightarrow \boldsymbol{\pi}^T G\boldsymbol{\alpha} = \boldsymbol{\beta}^T G\boldsymbol{\alpha} = \gamma + 1.$$

**Corollary 5.5.1**: *The parametric representation for G given by Theorem 5.5 is unique in that if G is any g-inverse of I – P then*

$$G(\boldsymbol{\alpha},\boldsymbol{\beta},\gamma) = [I-P+\boldsymbol{\alpha}\boldsymbol{\beta}^T]^{-1} + \gamma\boldsymbol{e}\boldsymbol{\pi}^T,$$

*where* $\boldsymbol{\alpha} = [I-(I-P)G]\boldsymbol{e}$, $\boldsymbol{\beta}^T = \boldsymbol{\pi}^T[I-G(I-P)]$, $\gamma+1 = \boldsymbol{\beta}^T G\boldsymbol{\alpha} = \boldsymbol{\beta}^T G\boldsymbol{e} = \boldsymbol{\pi}^T G\boldsymbol{\alpha}$.

The proof follows by noting that $\boldsymbol{\alpha} = A\boldsymbol{e}$, $\boldsymbol{\beta}^T = \boldsymbol{\pi}^T B$, as given by equation (5.16).

An important application (exercise – see (Hunter, 1990)) is that if $G$ has the representation $G(\boldsymbol{\alpha}, \boldsymbol{\beta}, \gamma)$:

$$G \in A\{1, 2\} \Leftrightarrow \gamma = -1, \quad (5.17)$$

$$G \in A\{1, 3\} \Leftrightarrow \boldsymbol{\alpha} = \boldsymbol{\pi}/\boldsymbol{\pi}^T\boldsymbol{\pi}, \quad (5.18)$$

$$G \in A\{1, 4\} \Leftrightarrow \boldsymbol{\beta} = \boldsymbol{e}/\boldsymbol{e}^T\boldsymbol{e} = \boldsymbol{e}/m, \quad (5.19)$$

$$G \in A\{1, 5\} \Leftrightarrow \boldsymbol{\alpha} = \boldsymbol{e}, \boldsymbol{\beta} = \boldsymbol{\pi}, \quad (5.20)$$

so that we have a simple procedure for classifying the g-inverses of $I - P$. As we shall see later, knowledge of the parameters $\boldsymbol{\alpha}$, $\boldsymbol{\beta}$, $\gamma$ can also provide useful information in determining the stationary distribution and moments of the first passage time distributions.

Special cases of generalized inverses of $I - P$ include the following:

(a) Kemeney and Snell's *fundamental matrix* of finite irreducible Markov chains,

$Z = [I - P + \boldsymbol{e}\boldsymbol{\pi}^T]^{-1} = G(\boldsymbol{e}, \boldsymbol{\pi}, 0)$, a (1, 5) g-inverse with $\gamma = 0$.

$Z = [I - P + \Pi]^{-1}$ where $\Pi = \boldsymbol{e}\boldsymbol{\pi}^T$ was introduced by (Kemeny and Snell, 1960).

Z was shown by was shown in (Hunter, 1969) to be a one condition g-inverse of $I - P$.

(b) The *Moore-Penrose* g-inverse of $I - P$ is the unique matrix $A^{(1,2,3,4)}$ satisfying conditions 1, 2, 3 and 4 is $G = G(\boldsymbol{\pi}/\boldsymbol{\pi}^T\boldsymbol{\pi}, \boldsymbol{e}/\boldsymbol{e}^T\boldsymbol{e}, -1)$ which can be expressed as

$$A^{(1,2,3,4)} = [I-P+\boldsymbol{\pi}\boldsymbol{e}^T]^{-1} - \frac{\boldsymbol{e}\boldsymbol{\pi}^T}{m\boldsymbol{\pi}^T\boldsymbol{\pi}}.$$

An equivalent form, $A^{(1,2,3,4)} = [I-P+\alpha\boldsymbol{\pi}\boldsymbol{e}^T]^{-1} - \alpha\boldsymbol{e}\boldsymbol{\pi}^T$, where $\alpha = 1/\sqrt{m\boldsymbol{\pi}^T\boldsymbol{\pi}}$ was originally derived by (Wachter, 1973), (Paige, Styan and Wachter, 1975). The equivalence of the expressions follows from equation (5.8).

(c) The *group inverse* of $I – P$ is the unique matrix $A^{(1,2,5)}$ satisfying conditions 1, 2 and 5 is $A^{\#} = G(\boldsymbol{e}, \boldsymbol{\pi}, -1) \equiv [I – P + \Pi]^{-1} - \Pi$. $A^{\#}$ was originally identified as the group inverse of $I – P$ in (Meyer, 1975).

(d) If $P = \begin{bmatrix} P_{11} & \alpha \\ \beta^T & p_{mm} \end{bmatrix}$ then $I – P$ has a g-inverse of the form

$$(I-P)^{-} = \begin{bmatrix} (I - P_{11})^{-1} & \boldsymbol{0} \\ \boldsymbol{0}^T & 0 \end{bmatrix} = [I - P + \boldsymbol{t}\boldsymbol{u}^T]^{-1} + \boldsymbol{e}\boldsymbol{f}^T,$$

where $\boldsymbol{u}^T = (\boldsymbol{0}^T, 1)$, $\boldsymbol{t}^T = (\boldsymbol{0}^T, 1)$, $\boldsymbol{f}^T = -(\boldsymbol{\beta}^T (I - P_{11})^{-1}, 1)$.

The partitioned form for $(I - P)^{-}$, is due to (Rhode, 1968).

## 6. Stationary Distributions

In many stochastic processes the determination of stationary distributions is an important problem as it leads, in the case of aperiodic irreducible processes, to knowledge of the limiting distribution.

We shall see later that the derivation of mean first passage times in Markov chains involves either the computation of a matrix inverse or a matrix g-inverse, so we consider only those techniques for solving the stationary distributions that use g-inverses. This will assist us later to consider the joint computation of the stationary distributions and mean first passage times with a minimal set of computations.

There are of course a variety of computational techniques that one can use to obtain a solution to the constrained singular system of linear equations $\boldsymbol{\pi}^T(I – P) = \boldsymbol{0}^T$, subject to the boundary condition $\boldsymbol{\pi}^T\boldsymbol{e} = 1$. We do not focus on such techniques in this presentation.

We consider three specific classes of procedures using generalized matrix inverses - one using $A = I – (I – P)G$, one using and $B = I – G (I – P)$, and one using simple particular forms of the g-inverse $G$ without computing either $A$ or $B$.

### 6.1 Procedures using $A$

From earlier, (equation (2.1)), we saw that the stationary probability vector, $\boldsymbol{\pi}^T = (\pi_1, \pi_2, ..., \pi_m)$, for a finite irreducible Markov chain satisfies the set of linear equations $\boldsymbol{\pi}^T(I - P) = \boldsymbol{0}^T$. This is an equation of type $XB = C$ with $X = \boldsymbol{\pi}^T$, $B = I – P$, $C = \boldsymbol{0}^T$ which must be solved subject to the restriction $\boldsymbol{\pi}^T\boldsymbol{e} = 1$. From equation (4.4) we saw that provided $CB^{-}B = C$ (which is obviously satisfied in this case), that $XB = C$ has the general solution $X = CB^{-} + W(I - BB^{-})$ where $W$ is an arbitrary matrix. Thus if $G$ is any generalised inverse of $I – P$, $\boldsymbol{\pi}^T = \boldsymbol{w}^T[I - (I - P)G] = \boldsymbol{w}^T A$, where $A = I – (I – P)G$ and $\boldsymbol{w}^T$ is chosen so that $\boldsymbol{\pi}^T\boldsymbol{e} = 1$. This leads to the following key result:

**Theorem 6.1**: (Hunter, 1982) *If G is any g-inverse of I – P, A ≡ I – (I – P)G and* $\mathbf{v}^T$ *is any vector such that* $\mathbf{v}^T A\mathbf{e} \neq 0$ *then*

$$\boldsymbol{\pi}^T = \frac{\mathbf{v}^T A}{\mathbf{v}^T A\mathbf{e}} \ . \tag{6.1}$$

*Furthermore* $A\mathbf{e} \neq 0$ *for all g-inverse of G so that it is always possible to find a suitable* $\mathbf{v}^T$.

**Proof**: From equation (5.14), $A = \boldsymbol{\alpha}\boldsymbol{\pi}^T$. so that $A\mathbf{e} = \boldsymbol{\alpha}\,(\neq 0)$ and for any $\mathbf{v}^T$, $\mathbf{v}^T A = \mathbf{v}^T\boldsymbol{\alpha}\boldsymbol{\pi}^T = \mathbf{v}^T(A\mathbf{e})\boldsymbol{\pi}^T = (\mathbf{v}^T A\mathbf{e})\boldsymbol{\pi}^T$, leading to equation (6.1).

From the proof of Theorem 5.5, $A = \boldsymbol{\alpha}\boldsymbol{\pi}^T$ where $\boldsymbol{\alpha}$ is a vector of the parameters specifying the g-inverse $G$. Thus it is clear that $\boldsymbol{\alpha}^T = (\alpha_1, \alpha_2, ..., \alpha_m) \neq \mathbf{0}$. This implies that there is at least one $i$ with $\alpha_i \neq 0$. Let $\mathbf{v}^T = (v_1, v_2, ..., v_m)$. Since $A\mathbf{e} = \boldsymbol{\alpha}$, $\mathbf{v}^T A\mathbf{e} = \mathbf{v}^T\boldsymbol{\alpha} = \sum_{k=1}^{m} v_k\alpha_k$. Suppose we take $\mathbf{v}^T = \mathbf{e}_i^T \equiv (0, ..., 0, 1, 0, ..., 0)$, with 1 in the $i$-th position, then $\mathbf{v}^T A\mathbf{e} = \alpha_i \neq 0$. Thus it is always possible to find a suitable $\mathbf{v}^T$ for Theorem 6.1. Knowledge of the conditions of the g-inverse usually leads to suitable choices of $\mathbf{v}^T$ that simplify $\mathbf{v}^T A\mathbf{e}$. Typical choices are $\mathbf{v}^T = \mathbf{e}^T$ or $\mathbf{e}_i^T$. Such choices often depend on the multi-condition g-inverse being used.

**Corollary 6.1.1**: (Hunter, 1992) *Let G be any g-inverse of I – P, and A= I – (I – P)G. Then*

$$\boldsymbol{\pi}^T = \frac{\mathbf{e}^T A^T A}{\mathbf{e}^T A^T A\mathbf{e}}. \tag{6.2}$$

**Proof**: Since $\boldsymbol{\alpha} = A\mathbf{e} \neq 0$ if $\mathbf{v} = \boldsymbol{\alpha}$ then $\mathbf{v}^T A\mathbf{e} = \boldsymbol{\alpha}^T\boldsymbol{\alpha} = \sum \alpha_i^2 > 0$ and the equation (6.2) follows from Theorem 6.1 with $\mathbf{v}^T = A\mathbf{e}$.

**Corollary 6.1.2**: *Let G be a (1,3) g-inverse of I – P, and A = I – (I – P)G. Then*

$$\boldsymbol{\pi}^T = \frac{\mathbf{e}^T A}{\mathbf{e}^T A\mathbf{e}} \ , \tag{6.3}$$

*and for any* $i = 1, 2, \ldots m$, if $\mathbf{e}_i^T$ *is the i-th elementary vector,*

$$\boldsymbol{\pi}^T = \frac{\mathbf{e}_i^T A}{\mathbf{e}_i^T A\mathbf{e}} \ . \tag{6.4}$$

**Proof:** If $G$ is a (1, 3) inverse, $\boldsymbol{\alpha} = \dfrac{\boldsymbol{\pi}}{\boldsymbol{\pi}^T\boldsymbol{\pi}}$ so that $A = \boldsymbol{\alpha}\boldsymbol{\pi}^T = \dfrac{\boldsymbol{\pi}\boldsymbol{\pi}^T}{\boldsymbol{\pi}^T\boldsymbol{\pi}} = A^T$. Thus $A$ is symmetric. But for all g-inverses, $A$ is idempotent i.e. $A^2 = A$ so that $A^T A = A^2 = A$ and equation (6.3) follows from Corollary 6.1.1. Furthermore, for all $i$, $\mathbf{e}_i^T A\mathbf{e} = \mathbf{e}_i^T\boldsymbol{\pi}\boldsymbol{\pi}^T\mathbf{e}/\boldsymbol{\pi}^T\boldsymbol{\pi} = \pi_i/\sum\pi_i^2 > 0$ and equation (6.4) follows from Theorem 6.1 with $\mathbf{v} = \mathbf{e}_i$.

(Decell and Odell, 1967) derived the expression (6.3) for a (1,3) g-inverse of $I - P$ while expression (6.4) is due to (Hunter, 1992).

**Corollary 6.1.3**: (Hunter, 1992) *Let G be a (1,5) g-inverse of I – P, and A = I – (I – P)G. Then*

$$\boldsymbol{\pi}^T = \frac{\boldsymbol{e}^T A}{\boldsymbol{e}^T \boldsymbol{e}} = \frac{\boldsymbol{e}^T A}{m} , \tag{6.5}$$

*and for any i* = 1, 2, …*m*, if $\boldsymbol{e}_i^T$ *is the i-th elementary vector,*

$$\boldsymbol{\pi}^T = \boldsymbol{e}_i^T A. \tag{6.6}$$

**Proof:** If *G* is a (1, 5) inverse, $\boldsymbol{\alpha} = \boldsymbol{e}$ so that $A = \boldsymbol{e}\boldsymbol{\pi}^T$. Consequently equation (6.5) follows from Theorem 6.1 by taking $\boldsymbol{v} = \boldsymbol{e}$ implying $\boldsymbol{v}^T A\boldsymbol{e} = \boldsymbol{e}^T\boldsymbol{e}\boldsymbol{\pi}^T\boldsymbol{e} = \boldsymbol{e}^T\boldsymbol{e} = m.$ Similarly, equation (6.6) follows by taking $\boldsymbol{v} = \boldsymbol{e}_i$ implying $\boldsymbol{v}^T A\boldsymbol{e} = \boldsymbol{e}_i^T\boldsymbol{e}\boldsymbol{\pi}^T\boldsymbol{e} = 1.$

(Meyer, 1975) established equation (6.5) under the assumption that *G* is a (1,2,5) inverse, i.e. the group inverse of $I - P$. Note however that the 2 condition is not necessary.

**Corollary 6.1.4**: (Rhode, 1968), (Meyer, 1975):

$$\text{If } P = \begin{bmatrix} P_{11} & \boldsymbol{\alpha} \\ \boldsymbol{\beta}^T & p_{mm} \end{bmatrix} \text{ then } \boldsymbol{\pi}^T = \frac{(\boldsymbol{\beta}^T(I - P_{11})^{-1}, 1)}{(\boldsymbol{\beta}^T(I - P_{11})^{-1}, 1)\mathbf{e}}. \tag{6.7}$$

**Proof**: With the g-inverse of $I - P$ of the form

$$G = \begin{bmatrix} (I - P_{11})^{-1} & \boldsymbol{0} \\ \boldsymbol{0}^T & 0 \end{bmatrix}, A = I - (I - P)G = \begin{bmatrix} 0 & \boldsymbol{0} \\ \boldsymbol{\beta}^T(I - P_{11})^{-1} & 1 \end{bmatrix}$$ using (6.1) with $\boldsymbol{v} = \boldsymbol{e}$, gives

$$\boldsymbol{\pi}^T = \frac{\boldsymbol{e}^T A}{\boldsymbol{e}^T A\boldsymbol{e}} = \frac{(\boldsymbol{\beta}^T(I - P_{11})^{-1}, 1)}{(\boldsymbol{\beta}^T(I - P_{11})^{-1}, 1)\mathbf{e}},$$ and equation (6.7) follows.

The following theorem gives a computational procedure that may usually only require the initial row of *A* computed.

**Theorem 6.2***:* (Hunter, 1992) *Let G be any g-inverse of* $I - P$. *Let* $A = I - (I - P)G \equiv [a_{ij}]$.
*Let r be the smallest integer i* ($1 \le i \le m$) *such that* $\sum_{k=1}^{m} a_{ik} \neq 0$, *then*

$$\pi_j = \frac{a_{rj}}{\sum_{k=1}^{m} a_{rk}}, \quad j = 1, 2, ..., m. \tag{6.8}$$

**Proof:** We have seen that $\alpha^T = (\alpha_1, \alpha_2, ..., \alpha_m) \neq 0$ so that there is at least one such *i* such that $\alpha_i \neq 0$. Now let $A = \boldsymbol{\alpha}\boldsymbol{\pi}^T = [a_{ij}]$ then $a_{ij} = \alpha_i \pi_j$, and since $\pi_j > 0$ for all *j*, $a_{ij}$ must be non-zero for at least one such *i*. Since $\boldsymbol{e}_i^T A = \boldsymbol{e}_i^T \boldsymbol{\alpha}\boldsymbol{\pi}^T = \alpha_i \boldsymbol{\pi}^T$ is the *i-th* row of *A*, we can always find at least one row of *A* that does not contain a non-zero element. Note that $\alpha_i = \sum_{i=1}^{m} a_{ik}$ implying that $\alpha_i$ is the sum of the elements in the *i*-th row of *A*. Furthermore, if there is at least one non-zero element in that row, all the elements in that row must be non-zero, since the rows of *A* are scaled versions of $\boldsymbol{\pi}^T$. In particular if $A = [a_{ij}]$ then there is at least one *i* such $a_{i1} \neq 0$ in which case $a_{ij} \neq 0$ for $j = 1, \ldots, m,$ hence that $a_{ij} = (\sum_{k=1}^{m} a_{ik})\pi_j$ This leads to equation (6.8) with $r$ taken as the smallest such $i$.

In applying Theorem 6.2 one needs to first find the g-inverse $G = [g_{ij}]$ and then compute $a_{11}$ ( $= 1 - g_{11} + \sum_{k=1}^{m} p_{1k} g_{k1}$). If $a_{11} \neq 0$ then the first row of $A$ will suffice to find the stationary probabilities. If not find $a_{21}$, $a_{31}$, … and stop at the first non-zero $a_{r1}$.

For some specific g-inverses we can establish that we need only find the first row of $A$. For example MATLAB uses the pseudo inverse routine pinv($I – P$), to generate the (1,2,3,4) g-inverse of $I – P$.

**Corollary 6.2.1**: (Hunter, 1992) *If G is a* (1, 3) *or* (1, 5) *g-inverse of* $I – P$, *and if* $A = I - (I - P)G \equiv [a_{ij}]$ *then*

$$\pi_j = \frac{a_{1j}}{\sum_{k=1}^{m} a_{1k}}, \quad j = 1, 2, \ldots, m. \tag{6.9}$$

**Proof**: If $G$ satisfies condition 3, $\boldsymbol{\alpha} = \boldsymbol{\pi}/\boldsymbol{\pi}^T\boldsymbol{\pi}$, in which case $\alpha_1 \neq 0$. Similarly if $G$ satisfies condition 5, $\boldsymbol{\alpha} = \boldsymbol{e}$ in which case $\alpha_1 = 1$. The non-zero form of $\alpha_1$ ensures $a_{11} \neq 0$.

G-inverse conditions 2 or 4 do not place any restrictions upon $\boldsymbol{\alpha}$ and consequently the non-zero nature of $a_{11}$ cannot be guaranteed in these situations.

## 6.2 Procedures using *B*

In certain cases the expression $B = I – G(I – P)$ can also be used to find an expression for $\boldsymbol{\pi}^T$.

**Theorem 6.3**: (Hunter, 1992) *Let G be any g-inverse of* $I – P$ *that is not a* (1, 2) *g-inverse. Let* $B = I – G(I – P)$ *and* $\boldsymbol{v}^T$ *any vector such that* $\boldsymbol{v}^T\boldsymbol{e} \neq 0$. *Then*

$$\boldsymbol{\pi}^T = \frac{\boldsymbol{v}^T BG}{\boldsymbol{v}^T BG\boldsymbol{e}}. \tag{6.10}$$

**Proof**: From equation (5.15), $B = \boldsymbol{e\beta}^T$ and $(\gamma+1)\boldsymbol{\pi}^T = \boldsymbol{\beta}^T G$. Thus, if $\boldsymbol{v}^T$ is any vector, then $\boldsymbol{v}^T BG = \boldsymbol{v}^T\boldsymbol{e\beta}^T G = (\boldsymbol{v}^T\boldsymbol{e})(\gamma+1)\boldsymbol{\pi}^T$ and $\boldsymbol{v}^T BG\boldsymbol{e} = (\boldsymbol{v}^T\boldsymbol{e})(\gamma+1)$ leading to the conclusion of the theorem since $\boldsymbol{v}^T\boldsymbol{e} \neq 0$ and $\gamma \neq -1$.

**Corollary 6.3.1**: (Hunter, 1992) *Let G be any g-inverse of* $I – P$, *and* $B = I – G(I – P)$.
*For all G, except a* (1, 2) *g-inverse of* $I – P$,

$$\boldsymbol{\pi}^T = \frac{\boldsymbol{e}^T BG}{\boldsymbol{e}^T BG\boldsymbol{e}} \tag{6.11}$$

*and, for any* $i = 1, 2, \ldots, m$,

$$\boldsymbol{\pi}^T = \frac{\boldsymbol{e}_i^T BG}{\boldsymbol{e}_i^T BG\boldsymbol{e}}. \tag{6.12}$$

**Proof**: From the proof of Theorem 6.3 it is easily seen that $\boldsymbol{v} = \boldsymbol{e}$ and $\boldsymbol{v} = \boldsymbol{e}_i$ are suitable choices since $\boldsymbol{e}^T\boldsymbol{e} = m \neq 0$ and $\boldsymbol{e}^T\boldsymbol{e}_i = 1$.

**Theorem 6.4**: (Hunter, 1992) *Let G be a* (1, 5) *g-inverse of I – P, and B = I – G* (*I – P*). *Then for any i* = 1, 2, …, *m,*

$$\boldsymbol{\pi}^T = \boldsymbol{e}_i^T B. \tag{6.13}$$

**Proof:** If $G$ is a (1, 5) inverse of $I - P$, $\boldsymbol{\pi}^T = \boldsymbol{\beta}^T = \boldsymbol{e}_i^T B$ for any $i$ = 1, 2, ..., $m$ and the result follows.

### 6.3 Procedures using *G*

If $G$ is of special structure one can often find an expression for $\boldsymbol{\pi}^T$ in terms of $G$ alone, without computing either $A$ or $B$. The added computation of $A$ or $B$ following the derivation of a g-inverse $G$ is typically unnecessary, especially when additional special properties of $G$ are given. For example:

**Theorem 6.5**: (Hunter, 1992) *If G is a* (1, 4) *g-inverse of I – P,* $\boldsymbol{\pi}^T = \dfrac{\boldsymbol{e}^T G}{\boldsymbol{e}^T G\boldsymbol{e}}$. (6.14)

**Proof:** If $G$ is a (1, 4) inverse of $I - P$, $\boldsymbol{\beta} = \boldsymbol{e}/\boldsymbol{e}^T\boldsymbol{e}$. Since $B = \boldsymbol{e}\boldsymbol{\beta}^T$, substitution and simplification of equation (6.11) yields the given expression.

If $G$ has been determined from a computer package that specifies the nature of the g-inverse (e.g. the Moore-Penrose g-inverse) then the structure of the parameters $\boldsymbol{\alpha}$, $\boldsymbol{\beta}$ and $\gamma$ are known from equations (5.17) – (5.20). Alternatively, if $\boldsymbol{\alpha}$, $\boldsymbol{\beta}$ and $\gamma$ have been calculated using Theorem 5.5, these parameters often lead directly to an expression for $\boldsymbol{\pi}^T$.
If $G$ is a (1, 2) inverse, $\gamma = -1$. This parameter provides no information about $\boldsymbol{\pi}^T$.
If $G$ is (1, 3) inverse, $\boldsymbol{\alpha} = \boldsymbol{\pi}/\boldsymbol{\pi}^T\boldsymbol{\pi}$ so that $\boldsymbol{\pi} = \boldsymbol{\alpha}/\boldsymbol{\alpha}^T\boldsymbol{\alpha}$, consistent with Corollary 6.1.2.
If $G$ is a (1, 4) inverse, $\boldsymbol{\beta} = \boldsymbol{e}/\boldsymbol{e}^T\boldsymbol{e}$, consistent with Theorem 6.5.
If $G$ is a (1, 5) inverse, $\boldsymbol{\pi} = \boldsymbol{\beta}$, consistent with Theorem 6.4.

Perhaps the simplest general procedure for determining the stationary distribution of a Markov chain, using any generalized inverse, is given by the following Theorem 6.6 that uses equation (5.2) of Theorem 5.1. Rather than classifying $G$ as a specific "multi-condition" g-inverse, we now focus on special class of g-inverses which are matrix inverses of the simple form $[I - P + \boldsymbol{t}\boldsymbol{u}^T]^{-1}$, where $\boldsymbol{t}$ and $\boldsymbol{u}^T$ are simple forms, selected to ensure that the inverse exists with $\boldsymbol{\pi}^T\boldsymbol{t} \neq 0$ and $\boldsymbol{u}^T\boldsymbol{e} \neq 0$. A general result for deriving an expression for $\boldsymbol{\pi}^T$ using such a g-inverse is the following.

**Theorem 6.6**: (Paige, Styan and Wachter, 1973), (Kemeny, 1981), (Hunter, 1982)
*If* $G = [I - P + \boldsymbol{t}\boldsymbol{u}^T]^{-1}$ *where* $\boldsymbol{u}$ *and* $\boldsymbol{t}$ *are any vectors such that* $\boldsymbol{\pi}^T\boldsymbol{t} \neq 0$ *and* $\boldsymbol{u}^T\boldsymbol{e} \neq 0$, *then*

$$\boldsymbol{\pi}^T = \frac{\boldsymbol{u}^T G}{\boldsymbol{u}^T G\boldsymbol{e}}. \tag{6.15}$$

*Hence, if* $G = [g_{ij}]$ and $\boldsymbol{u}^T = (u_1, u_2, \ldots, u_m)$,

$$\pi_j = \frac{\sum_{k=1}^{m} u_k g_{kj}}{\sum_{r=1}^{m} u_r \sum_{s=1}^{m} g_{rs}} = \frac{\sum_{k=1}^{m} u_k g_{kj}}{\sum_{r=1}^{m} u_r g_{r.}}, \; j = 1, 2, ..., m. \quad (6.16)$$

**Proof:** Using equation (5.2) it is easily seen that $\boldsymbol{u}^T[I - P + \boldsymbol{tu}^T]^{-1}\boldsymbol{e} = \boldsymbol{\pi}^T\boldsymbol{e}/\boldsymbol{\pi}^T\boldsymbol{t} = 1/\boldsymbol{\pi}^T\boldsymbol{t}$ and equation (6.15) follows. The elemental expression (6.16) follows from equation (6.15).

The form for $\boldsymbol{\pi}^T$ above has the added simplification that we need only determine $G$ (and not $A$ or $B$ as in Theorems 6.1 and 6.2 and their corollaries.) While it will be necessary to evaluate the inverse of the matrix $I - P + \boldsymbol{tu}^T$ this may either be the inverse of a matrix which has a simple special structure or the inverse itself may be one that has a simple structure. Further, we also wish to use this inverse to assist in the determination of the mean first passage times (see Section 7).

We consider special choices of $\boldsymbol{t}$ and $\boldsymbol{u}$ based either upon the simple elementary vectors $\boldsymbol{e}_i$, the unit vector $\boldsymbol{e}$, the rows and/or columns of the transition matrix $P$, and in one case a combination of such elements. Let $\boldsymbol{p}_a^{(c)} \equiv P\boldsymbol{e}_a$ denote the *a-th* column of $P$ and $\boldsymbol{p}_b^{(r)T} \equiv \boldsymbol{e}_b^T P$ denote the *b-th* row of $P$.

Table 1, (Hunter, 2007a), below lists of a variety of special g-inverses with their specific parameters.

Table 1: Special g-inverses

| Identifier | g-inverse | Parameters | | |
|---|---|---|---|---|
| | $[I - P + \boldsymbol{tu}^T]^{-1}$ | $\boldsymbol{\alpha}$ | $\boldsymbol{\beta}^T$ | $\gamma$ |
| $G_{ee}$ | $[I - P + \boldsymbol{ee}^T]^{-1}$ | $\boldsymbol{e}$ | $\boldsymbol{e}^T/m$ | $(1/m) - 1$ |
| $G_{eb}^{(r)}$ | $[I - P + \boldsymbol{e}\,\boldsymbol{p}_b^{(r)T}]^{-1}$ | $\boldsymbol{e}$ | $\boldsymbol{p}_b^{(r)T}$ | 0 |
| $G_{eb}$ | $[I - P + \boldsymbol{e}\,\boldsymbol{e}_b^T]^{-1}$ | $\boldsymbol{e}$ | $\boldsymbol{e}_b^T$ | 0 |
| $G_{ae}^{(c)}$ | $[I - P + \boldsymbol{p}_a^{(c)}\boldsymbol{e}^T]^{-1}$ | $\boldsymbol{p}_a^{(c)}/\pi_a$ | $\boldsymbol{e}^T/m$ | $(1/m\pi_a) - 1$ |
| $G_{ab}^{(c,r)}$ | $[I - P + \boldsymbol{p}_a^{(c)}\boldsymbol{p}_b^{(r)T}]^{-1}$ | $\boldsymbol{p}_a^{(c)}/\pi_a$ | $\boldsymbol{p}_b^{(r)T}$ | $(1/\pi_a) - 1$ |
| $G_{ab}^{(c)}$ | $[I - P + \boldsymbol{p}_a^{(c)}\boldsymbol{e}_b^T]^{-1}$ | $\boldsymbol{p}_a^{(c)}/\pi_a$ | $\boldsymbol{e}_b^T$ | $(1/\pi_a) - 1$ |
| $G_{ae}$ | $[I - P + \boldsymbol{e}_a\boldsymbol{e}^T]^{-1}$ | $\boldsymbol{e}_a/\pi_a$ | $\boldsymbol{e}^T/m$ | $(1/m\pi_a) - 1$ |
| $G_{ab}^{(r)}$ | $[I - P + \boldsymbol{e}_a\boldsymbol{p}_b^{(r)T}]^{-1}$ | $\boldsymbol{e}_a/\pi_a$ | $\boldsymbol{p}_b^{(r)T}$ | $(1/\pi_a) - 1$ |
| $G_{ab}$ | $[I - P + \boldsymbol{e}_a\boldsymbol{e}_b^T]^{-1}$ | $\boldsymbol{e}_a/\pi_a$ | $\boldsymbol{e}_b^T$ | $(1/\pi_a) - 1$ |
| $G_{tb}^{(c)}$ | $[I - P + \boldsymbol{t}_b\boldsymbol{e}_b^T]^{-1}$ $(\boldsymbol{t}_b \equiv \boldsymbol{e} - \boldsymbol{e}_b + \boldsymbol{p}_b^{(c)})$ | $\boldsymbol{t}_b$ | $\boldsymbol{e}_b^T$ | 0 |

All these results follow from the observation that if $G = [I - P + \boldsymbol{tu}^T]^{-1}$ then, from Corollary 5.5.1 and the proof of Theorem 5.1, the parameters are given by $\boldsymbol{\alpha} = \boldsymbol{t}/\boldsymbol{\pi}^T\boldsymbol{t}$, $\boldsymbol{\beta}^T = \boldsymbol{u}^T/\boldsymbol{u}^T\boldsymbol{e}$ and $\gamma + 1 = 1/\{(\boldsymbol{\pi}^T\boldsymbol{t})(\boldsymbol{u}^T\boldsymbol{e})\}$.

$G_{tb}^{(c)}$ is included in Table 1 as the update $\boldsymbol{t}_b\boldsymbol{e}_b^T$ replaces the *b-th* column of $I - P$ by $\boldsymbol{e}$. See (Paige, Styan and Wachter, 1975).

The special structure of the g-inverses given in Table 1 leads, in many cases, to very simple forms for the stationary probabilities.

In applying Theorem 6.6, observe that $\boldsymbol{\pi}^T = \boldsymbol{u}^T G$ if and only if $\boldsymbol{u}^T G\boldsymbol{e} = 1$ if and only if $\boldsymbol{\pi}^T\boldsymbol{t} = 1$. Simple sufficient conditions for $\boldsymbol{\pi}^T\boldsymbol{t} = 1$ are $\boldsymbol{t} = \boldsymbol{e}$ or $\boldsymbol{t} = \boldsymbol{\alpha}$ (cf. (5.14)). (This later condition is of use only if $\boldsymbol{\alpha}$ does not explicitly involve any of the stationary probabilities, as for $G_{tb}^{(c)}$ )

**Corollary 6.6.1**: *If* $G = [I - P + \boldsymbol{e}\boldsymbol{u}^T]^{-1}$ *where* $\boldsymbol{u}^T\boldsymbol{e} \neq 0$,

$$\boldsymbol{\pi}^T = \boldsymbol{u}^T G. \qquad (6.17)$$

*and hence if* $\boldsymbol{u}^T = (u_1, u_2, \ldots, u_m)$ *and* $G = [g_{ij}]$ *then*

$$\pi_j = \sum_{k=1}^{m} u_k g_{kj},\ j = 1, 2, ..., m. \qquad (6.18)$$

In particular, we have the following special cases:

(a) If $\boldsymbol{u}^T = \boldsymbol{e}^T$ then $G \equiv G_{ee} = [I - P + \boldsymbol{e}\boldsymbol{e}^T]^{-1} = [g_{ij}]$ and

$$\pi_j = \sum_{k=1}^{m} g_{kj} \equiv g_{\bullet j}. \qquad (6.19)$$

(b) If $\boldsymbol{u}^T = \boldsymbol{p}_b^{(r)T}$ then $G \equiv G_{eb}^{(r)} = [I - P + \boldsymbol{e}\boldsymbol{p}_b^{(r)T}]^{-1} = [g_{ij}]$ and

$$\pi_j = \sum_{k=1}^{m} p_{bk} g_{kj}. \qquad (6.20)$$

(c) If $\boldsymbol{u}^T = \boldsymbol{e}_b^T$ then $G \equiv G_{eb} = [I - P + \boldsymbol{e}\boldsymbol{e}_b^T]^{-1} = [g_{ij}]$ and

$$\pi_j = g_{bj}. \qquad (6.21)$$

(Paige, Styan and Wachter, 1973) recommended using the matrix expression leading to equation (6.19), $\boldsymbol{\pi}^T[I - P + \boldsymbol{e}\boldsymbol{u}^T] = \boldsymbol{u}^T$ with $\boldsymbol{u}^T = \boldsymbol{p}_b^{(r)T} \equiv \boldsymbol{e}_b^T P$ for some *b*. (See also Section 7).

**Corollary 6.6.2**: *If* $G = [I - P + \boldsymbol{t}\boldsymbol{e}^T]^{-1}$ *where* $\boldsymbol{\pi}^T\boldsymbol{t} \neq 0$,

$$\boldsymbol{\pi}^T = \frac{\boldsymbol{e}^T G}{\boldsymbol{e}^T G\boldsymbol{e}}, \qquad (6.22)$$

*and hence, if* $G = [g_{ij}]$, *then*

$$\pi_j = \frac{\sum_{k=1}^{m} g_{kj}}{\sum_{r=1}^{m}\sum_{s=1}^{m} g_{rs}} = \frac{g_{\bullet j}}{g_{\bullet\bullet}},\ j = 1, 2, ..., m. \qquad (6.23)$$

In particular, equation (6.23) holds for $G = G_{ae}^{(c)}$, $G_{ee}$ and $G_{ae}$.

In the special case of $G_{ee}$, using equations (5.2) or (5.3), it follows that $g_{\bullet\bullet} = 1$, and equation (6.23) reduces to (6.19).

**Corollary 6.6.3**: *If* $G = [I - P + \boldsymbol{t}\boldsymbol{e}_b^T]^{-1}$ *where* $\boldsymbol{\pi}^T\boldsymbol{t} \neq 0$,

$$\boldsymbol{\pi}^T = \frac{\boldsymbol{e}_b^T G}{\boldsymbol{e}_b^T G\boldsymbol{e}}, \qquad (6.24)$$

*and hence, if* $G = [g_{ij}]$, *then*

$$\pi_j = \frac{g_{bj}}{\sum_{s=1}^{m} g_{bs}} = \frac{g_{bj}}{g_{b\bullet}},\ j = 1, 2, ..., m. \qquad (6.25)$$

In particular, equations (6.24) hold for $G = G_{ab}^{(c)}$, $G_{ab}$, $G_{eb}$ and $G_{tb}^{(c)}$.

In the special cases of $G_{eb}$ and $G_{tb}^{(c)}$, $g_{b\bullet} = 1$ and equation (6.24) reduces to equation (6.21).

**Corollary 6.6.4**: *If* $G = [I - P + \boldsymbol{t}\boldsymbol{p}_b^{(r)T}]^{-1}$ *where* $\boldsymbol{\pi}^T\boldsymbol{t} \neq 0$,

$$\boldsymbol{\pi}^T = \frac{\boldsymbol{p}_b^{(r)T} G}{\boldsymbol{p}_b^{(r)T} G\boldsymbol{e}}, \qquad (6.26)$$

*and hence, if* $G = [g_{ij}]$, *then*

$$\pi_j = \frac{\sum_{k=1}^{m} p_{bk} g_{kj}}{\sum_{i=1}^{m} \sum_{s=1}^{m} p_{bi} g_{is}},\ j = 1, 2, ..., m. \qquad (6.27)$$

In particular, results (6.26) hold for $G = G_{ab}^{(c,r)}$, $G_{ab}^{(r)}$ and $G_{eb}^{(r)}$.

In the special case of $G_{eb}^{(r)}$, the denominator of (6.27) is 1 and equation (6.27) reduces to (6.20).

Thus we have been able to find simple elemental expressions for the stationary probabilities using any of the g-inverses in Table 1. In the special cases of $G_{ee}$, $G_{eb}^{(r)}$, $G_{eb}$ and $G_{tb}^{(c)}$ the denominator of the expression given by equations (6.16) is always 1. (In each other case, observe that denominator of the expression $\boldsymbol{u}^T G\boldsymbol{e}$ is in fact $1/\pi_b$, with $\boldsymbol{u}^T G = \boldsymbol{\pi}^T / \pi_b$.)

In (Hunter, 2007a) the g-inverses of Table 1 are considered in more detail in order to highlight their structure and special properties that may provide either a computational check or a reduction in the number of computations required.

Let $\boldsymbol{g}_a^{(c)} = G\boldsymbol{e}_a$ denote the *a-th* column of $G$ and $\boldsymbol{g}_b^{(r)T} = \boldsymbol{e}_b^T G$ denote the *b-th* row of $G$. From the definition of $G = [I - P + \boldsymbol{t}\boldsymbol{u}^T]^{-1}$, pre- and post-multiplication $I - P + \boldsymbol{t}\boldsymbol{u}^T$ yields

$$G - PG + \boldsymbol{t}\,\boldsymbol{u}^T G = I, \qquad (6.28)$$

$$G - GP + G\boldsymbol{t}\,\boldsymbol{u}^T = I. \qquad (6.29)$$

Pre-multiplication by $\boldsymbol{\pi}^T$ and post-multiplication by $\boldsymbol{e}$ yields the expressions given by equations (5.2) and (5.3), i.e. $\boldsymbol{u}^T G = \boldsymbol{\pi}^T / \boldsymbol{\pi}^T\boldsymbol{t}$ and $G\boldsymbol{t} = \boldsymbol{e}/\boldsymbol{u}^T\boldsymbol{e}$.

Relationships between the rows, columns and elements of $G$ were considered in (Hunter, 2007a)

Let $\boldsymbol{g}_{rowsum} = G\boldsymbol{e} = \sum_{j=1}^{m} \boldsymbol{g}_j^{(c)} = [g_{1\bullet}, g_{2\bullet}, ..., g_{m\bullet}]^T$ denote the column vector of row sums of $G$ and $\boldsymbol{g}_{colsum}^T = \boldsymbol{e}^T G = \sum_{j=1}^{m} \boldsymbol{g}_j^{(r)T} = [g_{\bullet 1}, g_{\bullet 2}, ..., g_{\bullet m}]$ the row vector of column sums of $G$.

Table 2 below is given in (Hunter, 2007a).

Table 2: Row and column properties of g-inverses

| $G$ *g-inverse* | $\boldsymbol{g}_a^{(c)}$ *a-th column* | $\boldsymbol{g}_{colsum}^T$ *Column sum* | $\boldsymbol{g}_b^{(r)T}$ *b-th row* | $\boldsymbol{g}_{rowsum}$ *Row sum* | *Other properties* |
|---|---|---|---|---|---|
| $G_{ee}$ | | $\boldsymbol{\pi}^T$ | | $\boldsymbol{e}/m$ | |
| $G_{eb}^{(r)}$ | | | $\boldsymbol{e}_b^T$ | $\boldsymbol{e}$ | $\boldsymbol{p}_b^{(r)T} G = \boldsymbol{\pi}^T$ |
| $G_{eb}$ | | | $\boldsymbol{\pi}^T$ | $\boldsymbol{e}$ | |
| $G_{ae}^{(c)}$ | $\boldsymbol{e}_a$ | $\boldsymbol{\pi}^T/\pi_a$ | | | $G\boldsymbol{p}_a^{(c)} = \boldsymbol{e}/m$ |
| $G_{ab}^{(c,r)}$ | $\boldsymbol{e}_a + (1 - p_{ba})\boldsymbol{e}$ | | | | $\boldsymbol{p}_b^{(r)T} G = \boldsymbol{\pi}^T/\pi_a$ <br> $G\boldsymbol{p}_a^{(c)} = \boldsymbol{e}$ |
| $G_{aa}^{(c)}(a = b)$ | $\boldsymbol{e}_a$ | | $\boldsymbol{\pi}^T/\pi_a$ | | |
| $G_{ab}^{(c)}(a \neq b)$ | $\boldsymbol{e} + \boldsymbol{e}_a$ | | $\boldsymbol{\pi}^T/\pi_a$ | | |
| $G_{ae}$ | $\boldsymbol{e}/m$ | $\boldsymbol{\pi}^T/\pi_a$ | | | |
| $G_{aa}^{(r)}(a = b)$ | $\boldsymbol{e}$ | | $\boldsymbol{e}_a^T$ | | $\boldsymbol{p}_b^{(r)T} G = \boldsymbol{\pi}^T/\pi_b$ |
| $G_{ab}^{(r)}(a \neq b)$ | $\boldsymbol{e}$ | | $\boldsymbol{e}_b^T + \boldsymbol{\pi}^T/\pi_a$ | | $\boldsymbol{p}_b^{(r)T} G = \boldsymbol{\pi}^T/\pi_a$ |
| $G_{ab}$ | $\boldsymbol{e}$ | | $\boldsymbol{\pi}^T/\pi_a$ | | |
| $G_{tb}^{(c)}$ | | | $\boldsymbol{\pi}^T$ | $\boldsymbol{e}_b$ | $G\boldsymbol{t}_b = \boldsymbol{e}$ |

A key observation is that stationary distribution can be found in terms of just the elements of the *b-th* row of $G_{eb}$, $G_{ab}^{(c)}$, $G_{ab}^{(r)}(a \neq b)$, $G_{ab}$ and $G_{tb}^{(c)}$. This requires the determination of just $m$ elements of $G$. We exploit these particular matrices later.

If the entire g-inverse has been computed the stationary distribution can be found in terms of $\boldsymbol{g}_{colsum}^T$, the row vector of column sums, in the case of $G_{ee}$, $G_{ae}^{(c)}$ and $G_{ae}$. In each of these cases there are simple constraints on $\boldsymbol{g}_a^{(c)}$ and $\boldsymbol{g}_{rowsum}$, possibly reducing the number of computations required, or at least providing a computational check.

In the remaining cases of $G_{eb}^{(r)}, G_{ab}^{(c,r)}$ and $G_{ab}^{(r)}$, the additional computation of $\boldsymbol{p}_b^{(r)T} G$ is required to lead to an expression for the stationary probabilities.

In (Hunter, 2007a) results were derived and given highlighting the differences between $G_{aa}^{(c)}$, $G_{aa}^{(r)}$ and $G_{aa}$. In fact $G_{aa}^{(c)}$ and $G_{aa}^{(r)}$ differ only in the $a$-th row and $a$-th column.

## 7. Moments of the first passage times

Let $M = [m_{ij}]$ be the mean first passage time matrix of a finite irreducible Markov chain with transition matrix $P$.

We saw in Section 2.2 that $m_{ij}$, the mean first passage time from state $i$ to state $j$, satisfies equation (2.3), $m_{ij} = 1 + \sum_{k \neq j} p_{ik} m_{kj}$, and the matrix M satisfies the matrix equation (2.4), namely

$$(I - P)M = E - PD, \tag{7.1}$$

where $D = M_d = (\Pi_d)^{-1}$ with $\Pi = \boldsymbol{e\pi}^T$, $E = \boldsymbol{ee}^T = [1]$.

**Lemma 7.1.** *If X is an arbitrary square matrix, and Λ is a diagonal matrix,*

$$(XE)_d = (X\Pi)_d D, \quad (X\Lambda)_d = X_d \Lambda, \quad E\Pi_d = \Pi.$$

**Proof:** These results are well known and easily derived. See (Hunter, 1983).

We now utilise the generalized method presented in Section 4 for solving of equations of the form of (7.1).

**Theorem 7.2**: (Hunter, 1982). *If G is any g-inverse of I – P, then*

$$M = [G\Pi - E(G\Pi)_d + I - G + EG_d]D. \tag{7.2}$$

**Proof:** Equation (2.1) is of the form $AX = C$, where $A = I - P$, $X = M$ and $C$ is known. From equation (4.5), consistent equations of this form can be solved using any g-inverse of $A$, $A^-$, with the general solution given by $X = A^-C + (I - A^-A)W$, where $W$ is an arbitrary matrix.

The consistency condition, $A^-AC = C$, i.e., $[I - (I - P)G](E - PD) = 0$ can be shown to be satisfied. (Exercise: Take any g-inverse of the form $G = [I - P + \boldsymbol{tu}^T]^{-1} + \boldsymbol{ef}^T + \boldsymbol{g\pi}^T$ and use the properties considered in the proof of Theorem 5.1).

Thus the general solution of equation (7.1), with $G$ is any g-inverse of $I - P$, is given by

$$M = G[E - PD] + [I - G(I - P)]W, \tag{7.3}$$

where $W$ is an arbitrary matrix. The arbitrariness of W can be eliminated by taking advantage of the knowledge of $D = M_d$.

We first simplify equation (7.3) by using the results given the proof of Theorem 5.5. For any g-inverse $G$ of $I - P$, from equations (5.11) and (5.15),

$$I - G(I - P) = B = \boldsymbol{e\beta}^T. \tag{7.4}$$

Thus, defining $\boldsymbol{\beta}^T W = \boldsymbol{w}^T$ and noting that $\boldsymbol{e}\boldsymbol{w}^T$ can be expressed as $E\Lambda$ say, where $\Lambda$ is a diagonal matrix, whose diagonal elements are those of the vector $\boldsymbol{w}^T$. Thus equation (7.3) can be expressed as

$$M = G[E - PD] + BW = GE - GPD + E\Lambda. \tag{7.5}$$

We now determine $\Lambda$ by taking the diagonal elements of equation (7.5) using Lemma 7.1, to obtain $D = (G\Pi)_d D - (GP)_d D + \Lambda$, implying that

$$\Lambda = [I - (G\Pi)_d + (GP)_d]D. \tag{7.6}$$

Since $E = \Pi D$, substitution of $\Lambda$, from equation (7.6) into equation (7.5) yields

$$M = [G\Pi - E(G\Pi)_d - GP + E(GP)_d + E]D. \tag{7.7}$$

Further simplification of equation (7.7) is possible. From equation (7.4) observe that $E(I - G + GP)_d = \boldsymbol{e}\boldsymbol{e}^T(\boldsymbol{e}\boldsymbol{\beta}^T)_d = \boldsymbol{e}\boldsymbol{\beta}^T = I - G + GP$, yielding, after further refinement, the required equation (7.2).

**Corollary 7.2.1**:

$$GE - E(G\Pi)_d D = M - [I - G + EG_d]D. \tag{7.8}$$

**Proof:** Result (7.8) follows from equation (7.1) by noting that $\Pi D = E$.

If $A = [a_{ij}]$ is a matrix, define $a_{i\bullet} = \sum_{j=1}^{m} a_{ij}$.

**Corollary 7.2.2**: (Hunter, 2007) *Under any of the following three equivalent conditions,*

(i) $G\boldsymbol{e} = g\boldsymbol{e}$, *g a constant,*

(ii) $GE - E(G\Pi)_d D = O$,

(iii) $G\Pi - E(G\Pi)_d = O$,

$$M = [I - G + EG_d]D. \tag{7.9}$$

**Proof:** Condition (i) $G\boldsymbol{e} = g\boldsymbol{e}$ implies $GE = G\boldsymbol{e}\boldsymbol{e}^T = g\boldsymbol{e}\boldsymbol{e}^T = gE$ and $E(G\Pi)_d D = E(G\boldsymbol{e}\boldsymbol{\pi}^T)_d D = E(g\boldsymbol{e}\boldsymbol{\pi}^T)_d D = gE\Pi_d D = gE$ leading to condition (ii). Since $\Pi D = E$ substitution in (ii) and post-multiplication by $D^{-1}$ leads to condition (iii). Further under condition (iii), $G\Pi = G\boldsymbol{e}\boldsymbol{\pi}^T = \boldsymbol{e}\boldsymbol{e}^T(G\Pi)_d$. Post-multiplication by $\boldsymbol{e}$, since $\boldsymbol{\pi}^T\boldsymbol{e} = 1$, yields $G\boldsymbol{e} = \boldsymbol{e}\boldsymbol{e}^T(G\Pi)_d\boldsymbol{e} = g\boldsymbol{e}$ where $g = \boldsymbol{e}^T(G\Pi)_d\boldsymbol{e}$, a constant $(= \sum_{k=1}^{m} g_{k\bullet}\pi_k)$. Thus (i) $\Rightarrow$ (ii) $\Rightarrow$ (iii) $\Rightarrow$ (i) and the conditions are equivalent. Result (7.9) follows from equation (7.2) under condition (ii). Note that condition (i) implies that $g_{k\bullet} = g$ for all $k$, and is equivalent to the generalised inverse $G$ having an eigenvalue $g$ with right eigenvector $\boldsymbol{e}$.

Elemental expressions for the $m_{ij}$ follow from Theorem 7.2 as follows.

**Corollary 7.2.3**: If $G = [g_{ij}]$ is any g-inverse of $I - P$,

$$m_{ij} = ([g_{jj} - g_{ij} + \delta_{ij}]/\pi_j) + (g_{i\bullet} - g_{j\bullet}), \text{ for all } i,j. \tag{7.10}$$

*Further, when* $G\boldsymbol{e} = g\boldsymbol{e}$,

$$m_{ij} = [g_{jj} - g_{ij} + \delta_{ij}]/\pi_j \quad \text{for all } i,j. \tag{7.11}$$

**Theorem 7.3:**

*(a)* *If* $Z = [I - P + \boldsymbol{e}\boldsymbol{\pi}^T]^{-1} = [z_{ij}]$ *then* $M = [m_{ij}] = [I - Z + EZ_d]D$, *and*

$$m_{ij} = \begin{cases} 1/\pi_j, & i = j; \\ (z_{jj} - z_{ij})/\pi_j, & i \neq j. \end{cases}$$

(b) *If* $A^{\#} = [I - P + \boldsymbol{e}\boldsymbol{\pi}^T]^{-1} - \boldsymbol{e}\boldsymbol{\pi}^T = [a^{\#}_{ij}]$ *then* $M = [m_{ij}] = [I - A^{\#} + EA_d^{\#}]D$ *and*

$$m_{ij} = \begin{cases} 1/\pi_j, & i = j; \\ (a^{\#}_{jj} - a^{\#}_{ij})/\pi_j, & i \neq j. \end{cases}$$

Since $\Pi = \boldsymbol{e}\boldsymbol{\pi}^T$, these are special cases of equation (7.9) with $G = Z$, Kemeny and Snell's fundamental matrix $Z = [I - P - \Pi]^{-1}$ (since $Z\boldsymbol{e} = \boldsymbol{e}$ and $g = 1$) as given initially in (Kemeny and Snell, 1960) and $G = A^{\#} = Z - \Pi$, Meyer's group inverse of $I - P$, (with $A^{\#}\boldsymbol{e} = \boldsymbol{0}$ and $g = 0$) as given by (Meyer, 1975). We identify them separately as they have to date been the primary methods of computing $M$. Note however that both (7.10) and (7.11) require the prior computation of the stationary probability vector in order to compute Z or $A^{\#}$.

The following joint computation procedure for $\pi_j$ and $m_{ij}$ was given in (Hunter, 1992) based upon Theorem 6.2 and Corollary 7.2.3 above. (The version below which appears in (Hunter, 2007a) corrects some minor errors given in the initial derivation.)

**Theorem 7.4:**

1. *Compute* $G = [g_{ij}]$, *be any g-inverse of* $I - P$.
2. *Compute sequentially rows* 1, 2, …$r$ ($\leq m$) *of* $A = I - (I - P)G \equiv [a_{ij}]$ *until* $\sum_{k=1}^{m} a_{rk}$, ($1 \leq r \leq m$) *is the first non-zero sum.*
3. *Compute* $\pi_j = \dfrac{a_{rj}}{\sum_{k=1}^{m} a_{rk}}, \quad j = 1, ..., m.$
4. *Compute* $m_{ij} = \begin{cases} \dfrac{\sum_{k=1}^{m} a_{rk}}{a_{rj}}, & i = j, \\ \dfrac{(g_{jj} - g_{ij})\sum_{k=1}^{m} a_{rk}}{a_{rj}} + \sum_{k=1}^{m}(g_{ik} - g_{jk}), & i \neq j. \end{cases}$

Note that the procedure contains the unnecessary additional computation of the elements of $A$. However, in general not all the elements of $A$ need to obtained and in most cases only the elements of the first row $a_{11}, a_{12}, ..., a_{1m}$ suffice. For such a scenario, all the mean first passage times can be expressed in terms of $m^2$ elements of the g-inverse $G = [g_{ij}]$ and the $m$ elements $\{a_{1k}\}$ ($k = 1, 2, ..., m$) of the first row of $A$.

We now consider using the special g-inverses given in Table 1 and 2 to find expressions for all the $\pi_j$ and the $m_{ij}$. The results are summarised in Table 3, based on (Hunter, 2007a).

Table 3: Joint computation of $\{\pi_j\}$ and $[m_{ij}]$ using special g-inverses

| g-inverse | $\pi_j$ | $m_{jj}$ | $m_{ij}\ (i \neq j)$ |
|---|---|---|---|
| $G_{ee}$ | $g_{\bullet j}$ | $1/g_{\bullet j}$ | $(g_{jj}-g_{ij})/g_{\bullet j}$ |
| $G_{eb}^{(r)}$ | $\sum_k p_{bk}g_{kj}$ | $1/\sum_k p_{bk}g_{kj}$ | $(g_{jj}-g_{ij})/\sum_k p_{bk}g_{kj}$ |
| $G_{eb}$ | $g_{bj}$ | $1/g_{bj}$ | $(g_{jj}-g_{ij})/g_{bj}$ |
| $G_{ae}^{(c)}, G_{ae}$ | $g_{\bullet j}/g_{\bullet\bullet}$ | $g_{\bullet\bullet}/g_{\bullet j}$ | $(g_{jj}-g_{ij})g_{\bullet\bullet}/g_{\bullet j}+(g_{i\bullet}-g_{j\bullet})$ |
| $G_{ab}^{(c,r)}, G_{ab}^{(r)}$ | $\dfrac{\sum_k p_{bk}g_{kj}}{\sum_i\sum_s p_{bi}g_{is}}$ | $\dfrac{\sum_i\sum_s p_{bi}g_{is}}{\sum_k p_{bk}g_{kj}}$ | $\dfrac{(g_{jj}-g_{ij})\sum_i p_{bi}g_{i\bullet}}{\sum_k p_{bk}g_{kj}}+(g_{i\bullet}-g_{j\bullet})$ |
| $G_{ab}^{(c)}, G_{ab}$ | $g_{bj}/g_{b\bullet}$ | $g_{b\bullet}/g_{bj}$ | $(g_{jj}-g_{ij})g_{b\bullet}/g_{bj}+(g_{i\bullet}-g_{j\bullet})$ |
| $G_{tb}^{(c)}$ | $g_{bj}$ | $1/g_{bj}$ | $(g_{jj}-g_{ij})/g_{bj}+(\delta_{bi}-\delta_{bj})$ |

The special case of $G_{eb}$ deserves highlighting.

**Theorem 7.5:** *If* $G_{eb} = [I - P + \boldsymbol{e}\,\boldsymbol{e}_b^T]^{-1} = [g_{ij}]$, *then*

$$\pi_j = g_{bj},\quad j = 1, 2, \ldots, m, \tag{7.12}$$

*and*

$$m_{ij} = \begin{cases} 1/g_{bj}, & i = j, \\ (g_{jj} - g_{ij})/g_{bj}, & i \neq j. \end{cases} \tag{7.13}$$

This is one of the simplest computational expressions for both the stationary probabilities and the mean first passage times for a finite irreducible Markov chain. These results do not appear to have been given any special attention in the literature.

If one wished to find a computationally efficient algorithm for finding $\pi_j$ based upon $G_{eb}$ then an alternative procedure would be to solve the equations $\boldsymbol{\pi}^T(I - P + \boldsymbol{e}\boldsymbol{e}_b^T) = \boldsymbol{e}_b^T$ directly. This reduces the problem to finding an efficient package for solving this system of linear equations. (Paige, Styan and Wachter, 1975) recommended solving for $\boldsymbol{\pi}$ using $\boldsymbol{\pi}^T(I - P + \boldsymbol{e}\boldsymbol{u}^T) = \boldsymbol{u}^T$ with $\boldsymbol{u}^T = \boldsymbol{e}_j^T P = \boldsymbol{p}_j^{(r)T}$, using Gaussian elimination with pivoting. Their other suggested choices included $\boldsymbol{u}^T = \boldsymbol{e}_j^T$, the recommended algorithm above.

It is interesting to observe that the particular matrix inverse we suggest for favourable consideration has been proposed in the past as the basis for a computational procedure for solving for the stationary probabilities. Since techniques for finding the $m_{ij}$ typically require the computation of a matrix inverse, $G_{eb}$ also appears to be a suitable candidate for this joint computation. (See also later, Theorem 7.9 for $M^{(2)}$.)

If $G = G_{tb}^{(c)} = [g_{ij}]$ then $m_{ij} = \dfrac{g_{jj} - g_{ij} + \delta_{ij}}{g_{bj}} + \delta_{bi} - \delta_{bj}$. These elemental expressions for $M$, appear in (Hunter, 1983) and, in the case $b = m$, in (Meyer, 1978).

In deriving the mean first passage times one is in effect solving the set of equations (2.3). If in this set of equations if we hold $j$ fixed, ($j = 1, 2, ..., m$) and let $\boldsymbol{m}_j^T = (m_{1j}, m_{2j}, ..., m_{mj})$ then equation (2.3) yields, (Hunter, 1983), $\boldsymbol{m}_j = [I - P + \boldsymbol{p}_j^{(c)}\boldsymbol{e}_j]^{-1}\boldsymbol{e} = G_{jj}^{(c)}\boldsymbol{e}$.

**Theorem 7.6:** *For fixed j,* $1 \le i \le m$,

(a) *If* $G_{jj}^{(c)} = [g_{rs}]$ *then* $m_{ij} = \boldsymbol{e}_i^T G_{jj}^{(c)}\boldsymbol{e} = g_{i\bullet}$.

(*b*)*If* $G_{jj}^{(r)} = [g_{rs}]$ *then* $m_{ij} = \boldsymbol{e}_i^T G_{jj}^{(r)}\boldsymbol{e} + \dfrac{\delta_{ij}}{\pi_j} - 1 = g_{i\bullet} + \dfrac{\delta_{ij}}{\pi_j} - 1 = \begin{cases} \sum_{k=1}^{m} p_{jk} g_{k\bullet}, & i = j, \\ g_{i\bullet} - 1, & i \ne j. \end{cases}$

(*c*) *If* $G_{jj} = [g_{rs}]$ *then* $m_{ij} = \boldsymbol{e}_i^T G_{jj}\boldsymbol{e} + \dfrac{\delta_{ij} - 1}{\pi_j} = \begin{cases} g_{j\bullet} & i = j, \\ g_{i\bullet} - g_{j\bullet} & i \ne j. \end{cases}$

**Proof**: (Exercise – See (Hunter, 2007a).

The utilisation of special matrix inverses as g-inverses in the joint computation of stationary distributions and mean first passage times leads to a significant simplification in that at most a single matrix inverse often needs to be computed and typically this involves a row or column sum with a very simple form, further reducing the necessary computations. While no computational examples have been included in this presentation, a variety of new procedures have been given that warrant further examination from a computational efficiency perspective.

We parallel the development, given above for $M$, for $M^{(2)}$.

**Theorem 7.7**: $M^{(2)}$ *satisfies the matrix equation*

$$(I - P)M^{(2)} = E + 2P(M - M_d) - PM_d^{(2)}. \tag{7.14}$$

**Proof:** This result is well known. See (Hunter, 1982), (Hunter, 1983).

**Corollary 7.7.1**:

$$M_d^{(2)} = 2D(\Pi M)_d - D. \tag{7.15}$$

**Proof:** Pre-multiplication of both sides of equation (7.14) with $\Pi$, noting that $M_d = D$, $\Pi = \Pi P$, and $\Pi E = E$, yields $\Pi M_d^{(2)} = E + 2\Pi(M - D)$. Now take diagonal elements, using Lemma 7.1, to obtain $\Pi_d M_d^{(2)} = I + 2(\Pi M)_d - 2\Pi_d D$. Equation (7.15) follows since $\Pi_d^{-1} = D$.

**Corollary 7.7.2**: *If G is any g-inverse of* $I - P$,

$$M_d^{(2)} = D + 2D\{(I - \Pi)G(I - \Pi)\}_d D. \tag{7.16}$$

**Proof:** From equation (2.3), noting that $\Pi E = E$, it is easily seen that

$$\Pi M - M = [(I - \Pi)G(I - \Pi) + \Pi - I]D.$$

Taking the diagonal elements, and noting that $M_d = D$ and $\Pi_d D = I$, yields

$$(\Pi M)_d = \{(I - \Pi)G(I - \Pi)\}_d D + I. \tag{7.17}$$

Equation (7.16) now follows from equation (7.14).

**Corollary 7.7.3**: *If* $Ge = ge$,

$$M_d^{(2)} = D + 2D\{(I - \Pi)G\}_d D. \tag{7.18}$$

*In particular,*

$$M_d^{(2)} = D + 2DA_d^{\#}D, \tag{7.19}$$

$$= 2DZ_dD - D. \tag{7.20}$$

**Proof:** Equation (7.18) follows from equation (7.16) by observing that $G\Pi = Ge\boldsymbol{\pi}^T = ge\boldsymbol{\pi}^T = g\Pi$ and that $(I - \Pi)\,\Pi = 0$ since $\Pi^2 = \boldsymbol{e\pi}^T\boldsymbol{e\pi}^T = \boldsymbol{e\pi}^T = \Pi$.
Equation (7.19) follows directly from Equation (7.16) as it has been shown (Hunter, 1982) that for all g-inverses $G$ of $I - P$, $(I - \Pi)G(I - \Pi)$ is invariant and is in fact $A^{\#}$, the group inverse of $I - P$.
Equation (7.20) follows directly from Equation (7.19)) since $A^{\#} = Z - \Pi$.

Expression (7.16) is also given in (Hunter, 1990); expression (7.19) is given in (Meyer, 1975), and (Hunter, 1983) and expression (7.20) is given in (Kemeny and Snell, 1960).

The results given by the following theorem offer possible computational advantages.

**Theorem 7.8**: (Hunter, 2007b) *If $G$ is any g-inverse of $I - P$, then*

$$M^{(2)} = 2[GM - E(GM)_d] + [I - G + EG_d][M_d^{(2)} + D] - M, \tag{7.21}$$

$$= 2[GM - E(GM)_d] + 2[I - G + EG_d]D(\Pi M)_d - M. \tag{7.22}$$

**Proof:** Equation (7.14) is of the form $AX = C$, where $A = I - P$, $X = M^{(2)}$ and $C$ is known. Consistent equations of this form can be solved using any g-inverse of $A$, $A^-$, with the general solution given by $X = A^-B + (I - A^-A)U$, where $U$ is an arbitrary matrix (Hunter, 1982). (The consistency condition, $A^-AC = C$, can be shown to be satisfied.)
The general solution of equation (7.14), with $G$ is any g-inverse of $I - P$, is given by

$$M^{(2)} = G[E + 2P(M - M_d) - PM_d^{(2)}] + [I - G(I - P)]U \tag{7.23}$$

where $U$ is an arbitrary matrix. The arbitrariness of $U$ can be eliminated by taking advantage of the knowledge of $M_d^{(2)}$. We first simplify equation (7.23) by using equation (7.14) to show that

$$GP(M - M_d) = GPM - GPD = GM - GE.$$

Secondly, (as in the proof of Theorem 7.2), for any g-inverse $G$ of $I - P$, from equations (5.10) and (5.14),

$$I - G(I - P) = \boldsymbol{e\beta}^T. \tag{7.24}$$

Thus, defining $\boldsymbol{\beta}^T U = \boldsymbol{b}^T$ and noting that $\boldsymbol{eb}^T$ can be expressed as $EB$, say, where $B$ is a diagonal matrix, whose diagonal elements are those of the vector $\boldsymbol{b}^T$. Thus equation (7.23) can be expressed as

$$M^{(2)} = 2GM - GE - GPM_d^{(2)} + EB. \qquad (7.25)$$

We now determine $B$ by taking the diagonal elements of equation (7.25) and, using an appropriate expression for $M_d^{(2)}$, to obtain

$$B = M_d^{(2)} - 2(GM)_d + (G\Pi)_d D + (GP)_d M_d^{(2)}. \qquad (7.26)$$

Substitution of $B$, from equation (7.26) into equation (7.25) yields

$$M^{(2)} = 2GM - 2E(GM)_d - GE + E(G\Pi)_d D + [E + E(GP)_d - GP]M_d^{(2)}. \qquad (7.27)$$

Further simplification of equation (7.27) is possible. From equation (7.24) (as in the proof of Theorem 7.2) observe that $E(I - G + GP)_d = \boldsymbol{e}\boldsymbol{e}^T(\boldsymbol{e}\boldsymbol{\beta}^T)_d = \boldsymbol{e}\boldsymbol{\beta}^T = I - G + GP$, yielding,

$$M^{(2)} = 2GM - 2E(GM)_d - GE + E(G\Pi)_d D + [I - G + EG_d]M_d^{(2)}. \qquad (7.28)$$

We can now make use of Corollary 7.2.1 to obtain equation (7.21).
Equation (7.22) follows from Corollary 7.7.1.

**Corollary 7.8.1**: (Hunter, 2007b) *If* $G$ *is any g-inverse of* $I - P$*, such that* $G\boldsymbol{e} = g\boldsymbol{e}$*, then*

$$M^{(2)} = 2[GM - E(GM)_d] + MD^{-1}M_d^{(2)}. \qquad (7.29)$$

*and hence*

$$M^{(2)} = 2[ZM - E(ZM)_d] + M(2Z_dD - I), \qquad (7.30)$$

$$= 2[A^{\#}M - E(A^{\#}M)_d] + M(2A_d^{\#}D + I). \qquad (7.31)$$

**Proof.** When $G\boldsymbol{e} = g\boldsymbol{e}$, from Corollary 7.4.2, $GE - E(G\Pi)_d D = 0$ and $[I - G + EG_d] = MD^{-1}$. Equation (7.29) now follows immediately from equation (7.28). Equations (7.30) and (7.31) follow from equation (7.29) making use of equations (7.19)) and (7.20).

The original derivation of equation (7.30) was due to (Kemeny and Snell, 1960) but was derived by indirect methods. A proof using $Z$, along the lines given above, but not using arbitrary g-inverses, was given in (Hunter, 1983). (Meyer, 1975) gave an indirect proof using $A^{\#}$ leading to equation (7.31).

Equations (7.30) and (7.31) have, in the past, provided the standard computational methods for finding $M^{(2)}$. Theorem 7.8 has an added computational advantage in that any g-inverse of $I - P$ can be used. This is particularly important since the previously used g-inverses ($Z$ and $A^{\#}$) both required and involved expressions for the stationary distribution of the Markov chain.

The efficient computation of the moments of first passage times in a Markov chain has previously attracted interest, for example (Heyman and Reeves, 1989), and (Heyman and O'Leary, 1995). Techniques for solving linear equations occurring in Markov chains involving Gaussian elimination, algorithmic methods, state reduction and hybrid methods have been popular but in many instances they have been based upon the closed form solutions given by Equation (7.9) for the mean first passage times, and equations (7.30) and (7.31) for the second moments and variances, using either $Z$ or $A^{\#}$. In particular it was noted in (Heyman and Reeves, 1989) that there are potential accuracy problems in using such closed forms where it is pointed out that "the first thing that needs to be done is to compute $Z$. The majority of the work is to compute $\boldsymbol{\pi}$ and to do a matrix inversion." Further "There are three sources of numerical error. The first is the algorithm to compute $\boldsymbol{\pi}$. The second occurs in computing the inverse of $(I - P - \Pi)$; this matrix may have negative elements, and this can cause round-off errors when the inverse is evaluated. The third is the matrix multiplication in equation (7.9); the matrix multiplying $D$ may have negative elements. Now we consider the additional work to compute $M^{(2)}$, and the additional numerical errors that might occur. There are three matrix multiplications that are required, two of which involve at least one diagonal matrix. … In each of these multiplications there is a matrix with (possibly) negative elements, which may introduce round-off errors." In (Heyman and O'Leary, 1995) it is noted that "Deriving means and variances of first passage times from either the fundamental matrix $Z$ or the group generalised inverse $A^{\#}$ leads to a significant inaccuracy on the more difficult problems." The authors then conclude that "for this reason, it does not make sense to compute either the fundamental matrix or the group generalised inverse unless the individual elements of those matrices are of interest."

An alternative approach to exploring the moments of first passage times in a discrete time Markov chain is to consider them as the moments of a discrete-time phase type distribution that arise as times to absorption in an absorbing Markov chain. This approach, pioneered by (Neuts, 1981) has some useful computational procedures for fixed initial starting state $i$ and final state $j$ (as an absorbing state). Basically this involves the matrix inversion of a submatrix obtained from the transition matrix $P$ through the deletion of the $j$-th column and row, i.e. through $(I - P)_j$. This inverse can in fact be expressed as a special generalised inverse of $I - P$, (Hunter, 1983). Expressions for higher moments of the first passage times arise as factorial moments involving powers of the inverse of $(I - P)_j$ . We have not followed this line of attack in this paper as it typically focuses on the initial state $i$ (through the initial probability vector) and final state $j$, whereas the approach taken above leads to omnibus expressions for the first and second moments. The approach used in this paper enables comparison with the procedures pioneered in (Kemeny and Snell, 1960) and extended in (Meyer, 1975). The only closed form expressions utilised in the literature for finding $M$ and $M^{(2)}$, thus far, have been expressions involving $Z$ and $A^{\#}$.

Corollary 7.8.1 provides a much simpler form for the computation of $M^{(2)}$ if one is prepared to restrict attention to the class of g-inverses of $I - P$ that have the property $G\boldsymbol{e} = g\boldsymbol{e}$. Note that if $G = (I - P)_j = [I - P + \boldsymbol{t}_j\boldsymbol{e}_j^T]^{-1}$ with $\boldsymbol{t}_j = \boldsymbol{e} - (I - P)\boldsymbol{e}_j$, as used in the phase type distribution approach, then $G\boldsymbol{e} = \boldsymbol{e}_j$, see (Hunter, 1983) so that this restriction property is not satisfied. While $Z$ and $A^{\#}$ both satisfy this property we can take advantage of much simpler forms of such $g$-inverses. We explore some further consequences of this observation in Theorem 7.9 below.

Elemental expressions for $m_{ij}^{(2)}$ can be found from equation (7.21). Note that if $A = [a_{ij}]$ then $EA_d = [a_{jj}]$.

**Corollary 7.8.2**: *If $G = [g_{ij}]$ is any g-inverse of $I - P$ then*

$$m_{ij}^{(2)} = 2\sum_{k=1}^{m}(g_{ik} - g_{jk})m_{kj} - m_{ij} + (\delta_{ij} - g_{ij} + g_{jj})(m_{jj}^{(2)} + m_{jj}). \qquad (7.32)$$

*Further, if $G\boldsymbol{e} = g\boldsymbol{e}$, then*

$$m_{ij}^{(2)} = 2\sum_{k=1}^{m}(g_{ik} - g_{jk})m_{kj} + m_{ij}m_{jj}^{(2)} \big/ m_{jj}. \qquad (7.33)$$

Equations (7.32) and (7.33) do not provide an explicit expressions for the $m_{jj}^{(2)}$. These can be derived from equation (7.15).

**Corollary 7.8.3**:

$$m_{jj}^{(2)} + m_{jj} = 2m_{jj}\sum_{i=1}^{m}\pi_i m_{ij}. \qquad (7.34)$$

To utilise equation (7.41), simplified expressions of $\sum_{i=1}^{m} \pi_i m_{ij}$ are required.

**Corollary 7.8.4**: *If* $\boldsymbol{\alpha}^T \equiv (\alpha_1,...,\alpha_m)$ *where* $\alpha_j \equiv \sum_{i=1}^{m} \pi_i m_{ij}$, *then*

$$\boldsymbol{\alpha}^T = \boldsymbol{\pi}^T M = \boldsymbol{e}^T (\Pi M)_d. \tag{7.35}$$

*Further,* $\boldsymbol{\alpha} = (\Pi M)_d \boldsymbol{e}$ *and if G is any g-inverse of I – P, then*

$$\boldsymbol{\alpha} = (\boldsymbol{\pi}^T G\boldsymbol{e})\boldsymbol{e} - (GE)_d \boldsymbol{e} + \boldsymbol{e} - (\Pi G)_d D\boldsymbol{e} + G_d D\boldsymbol{e}. \tag{7.36}$$

*In particular if* $G\boldsymbol{e} = g\boldsymbol{e}$,

$$\boldsymbol{\alpha} = \boldsymbol{e} - (\Pi G)_d D\boldsymbol{e} + G_d D\boldsymbol{e}. \tag{7.37}$$

*If Z is the fundamental matrix and* $A^\#$ *is the group inverse,*

$$\boldsymbol{\alpha} = Z_d D\boldsymbol{e} = \boldsymbol{e} + A_d^\# D\boldsymbol{e}. \tag{7.38}$$

**Proof**: The first expression of equation (7.35) follows from the definition. Note that $\Pi M = \boldsymbol{e}\boldsymbol{\pi}^T M = \boldsymbol{e}\boldsymbol{\alpha}^T$, and consequently $\alpha_j = (\Pi M)_{jj} = (\Pi M)_{ij}$ for all *i, j.* This implies that $\boldsymbol{\alpha}^T = \boldsymbol{e}^T(\Pi M)_d$ and hence $\boldsymbol{\alpha} = (\Pi M)_d \boldsymbol{e}$.

From equation (7.8), $\boldsymbol{\alpha}^T = \boldsymbol{\pi}^T M = \boldsymbol{\pi}^T [G\Pi - E(G\Pi)_d + I - G + EG_d]D$. Simplification yields $\boldsymbol{\alpha}^T = (\boldsymbol{\pi}^T G\boldsymbol{e})\boldsymbol{e}^T - \boldsymbol{e}^T (GE)_d + \boldsymbol{e}^T - \boldsymbol{e}^T (\Pi G)_d D + \boldsymbol{e}^T G_d D$, since $\Pi D = E = \boldsymbol{e}\boldsymbol{e}^T$, $\boldsymbol{\pi}^T E = \boldsymbol{e}^T$, $\boldsymbol{\pi}^T D = \boldsymbol{e}^T$, $\boldsymbol{\pi}^T G = \boldsymbol{e}^T (\Pi G)_d$, and from equation (7.2), $(G\Pi)_d E = (GE)_d$. Equation (7.36) follows by noting that if $\Lambda$ is a diagonal matrix and $\boldsymbol{a}^T = \boldsymbol{e}^T\Lambda$ then $\boldsymbol{a} = \Lambda\boldsymbol{e}$.

An alternative derivation also follows from equation (7.17) by noting that $\boldsymbol{\alpha} = (\Pi M)_d \boldsymbol{e} = [(I - \Pi)G(I - \Pi)]_d D\boldsymbol{e} + \boldsymbol{e}$. The equivalence follows by noting that $(\Pi G \Pi)_d D = (\boldsymbol{\pi}^T G\boldsymbol{e})I$.

When $G\boldsymbol{e} = g\boldsymbol{e}$, $(\boldsymbol{\pi}^T G\boldsymbol{e})\boldsymbol{e} - (G\Pi)_d D\boldsymbol{e} = g(\boldsymbol{\pi}^T \boldsymbol{e})\boldsymbol{e} - g\Pi_d D\boldsymbol{e} = g\boldsymbol{e} - g\boldsymbol{e} = \boldsymbol{0}$, and equation (7.37) follows from equation (7.36). Equations (7.38) follow from equation (7.37) and the facts that $(\Pi Z)_d D = \Pi_d D = I$ and $A^\# = Z - \Pi$.

The transpose variant of equation (7.38), using *Z,* has been derived earlier - see (Hunter, 1983).

Elemental expressions for the $\alpha_j$ follow either from the expressions of Corollary 7.8.6 for $\boldsymbol{\alpha}^T$.

**Corollary 7.8.5**: *If* $G = [g_{ij}]$ *is any g-inverse of I – P, then*

$$\alpha_j = \sum_{i=1}^{m} \pi_i m_{ij} = 1 + \sum_{i=1}^{m} \pi_i g_{i\bullet} - g_{j\bullet} + \left(g_{jj} - \sum_{i=1}^{m} \pi_i g_{ij}\right)\Big/\pi_j.$$

*In particular if* $G\boldsymbol{e} = g\boldsymbol{e}$ *(i.e.* $g_{i\bullet} = g$ *for all i),* $\alpha_j = 1 + \left(g_{jj} - \sum_{i=1}^{m} \pi_i g_{ij}\right)\Big/\pi_j$.

*Further if* $Z = [z_{ij}]$ *and* $A^\# = [t_{ij}]$, $\alpha_j = \left(z_{jj}/\pi_j\right) = 1 + \left(t_{jj}/\pi_j\right)$.

We return to Theorem 7.5 and extend our analysis in the case when $G_{eb} = [I - P + \boldsymbol{e}\boldsymbol{e}_b^T]^{-1}$. We have already seen that this particular g-inverse has many desirable characteristics. In particular it was shown that if $G = G_{eb} = [g_{ij}]$ then (equations (7.11) and (7.11), respectively):

$$\pi_j = g_{bj},\ \ j = 1,2,...,m, \quad \text{and} \quad m_{ij} = (\delta_{ij} + g_{jj} - g_{ij})/g_{bj} = \begin{cases} 1/g_{bj}, & i = j, \\ (g_{jj} - g_{ij})/g_{bj}, & i \neq j. \end{cases}$$

We now use the same matrix $G = G_{eb}$ to also find expressions for $M_d^{(2)}$ and $M^{(2)}$.

**Theorem 7.9**: *If* $G = G_{eb} = [I - P + \boldsymbol{e}\boldsymbol{e}_b{}^T]^{-1} = [g_{ij}]$, *then*

$$m_{ij}^{(2)} = \begin{cases} m_{jj}[1 + 2m_{jj}(g_{jj} - g_{bj}^{(2)})], & i = j, \\ 2m_{jj}[g_{jj}^{(2)} - g_{ij}^{(2)} + m_{ij}(g_{jj} - g_{bj}^{(2)})] - m_{ij}, & i \neq j. \end{cases} \tag{7.39}$$

*Further*

$$var[T_{ij}] = \begin{cases} m_{jj}[1 - m_{jj} + 2m_{jj}(g_{jj} - g_{bj}^{(2)})], & i = j, \\ 2m_{jj}[g_{jj}^{(2)} - g_{ij}^{(2)} + m_{ij}(g_{jj} - g_{bj}^{(2)})] - m_{ij}(1 - m_{ij}), & i \neq j. \end{cases} \tag{7.40}$$

**Proof**: Firstly, from Corollary 7.7.2, $M = D - GD + EG_dD,$ where $D = [\delta_{ij}/\pi_j] = [\delta_{ij}/g_{bj}] = diag(1/g_{b1}, 1/g_{b2}, ..., 1/g_{bm})$. Thus $GD = [g_{ij}/g_{bj}], EG_dD = [g_{jj}/g_{bj}]$

Further $\Pi G = \boldsymbol{e}\boldsymbol{\pi}^T G = \boldsymbol{e}\boldsymbol{e}_{\boldsymbol{b}}^T G^2 = [\sum_{k=1}^m g_{bk}g_{kj}]$ so that $(\Pi G)_d = [\delta_{ij}\sum_{k=1}^m g_{bk}g_{kj}] = [\delta_{ij}g_{bj}^{(2)}] =$
$= diag(g_{b1}^{(2)}, g_{b2}^{(2)}, ..., g_{bm}^{(2)})$. From Corollary 7.7.3, $M_d^{(2)} = D + 2G_dD^2 - 2(\Pi G)_dD^2$ where $G_dD^2 = [\delta_{ij}g_{jj}/g_{bj}^2]$, $\Pi G = \boldsymbol{e}\boldsymbol{e}_{\boldsymbol{b}}^T G^2 = [g_{bj}^{(2)}]$ and hence $(\Pi G)_dD^2 = [\delta_{ij}\, g_{bj}^{(2)}/g_{bj}^2]$ implying that $m_{jj}^{(2)} = (1/g_{bj}) + 2(g_{jj}/g_{bj}^2) - 2(g_{bj}^{(2)}/g_{bj}^2) = [g_{bj} + 2(g_{jj} - g_{bj}^{(2)})]/g_{bj}^2$.
Thus, from Equation (7.11),

$$m_{jj}^{(2)} = m_{jj} + 2m_{jj}^2(g_{jj} - g_{bj}^{(2)}). \tag{7.41}$$

Note that since $G\boldsymbol{e} = \boldsymbol{e}$, $\sum_{k=1}^m g_{ik} = g_{i\bullet} = 1$, and thus

$$\begin{aligned} \sum_{k=1}^m (g_{ik} - g_{jk})m_{kj} &= [(g_{ij} - g_{jj})m_{jj}] + [\sum_{k \neq j}^m (g_{ik} - g_{jk})(g_{jj} - g_{kj})m_{jj}] \\ &= (g_{ij} - g_{jj})m_{jj} + \sum_{k=1}^m (g_{ik} - g_{jk})(g_{jj} - g_{kj})m_{jj} \\ &= m_{jj}[g_{ij} - g_{jj} + \sum_{k=1}^m (g_{ik}g_{jj} - g_{ik}g_{kj} - g_{jk}g_{jj} + g_{jk}g_{kj})] \\ &= m_{jj}[g_{ij} - g_{jj} + (\sum_{k=1}^m g_{ik})g_{jj} - g_{ij}^{(2)} - (\sum_{k=1}^m g_{jk})g_{jj} + g_{jj}^{(2)}] \\ &= m_{jj}[g_{ij} - g_{jj} + g_{i\bullet}g_{jj} - g_{ij}^{(2)} - g_{j\bullet}g_{jj} + g_{jj}^{(2)}] = m_{jj}[g_{ij} - g_{jj} + g_{jj}^{(2)} - g_{ij}^{(2)}]. \end{aligned}$$

Substituting in equation (7.33), for $i \neq j$, and using equation (7.11), yields

$$m_{ij}^{(2)} = 2m_{jj}[g_{jj}^{(2)} - g_{ij}^{(2)}] + m_{ij}\left\{\frac{m_{jj}^{(2)}}{m_{jj}} - 2\right\}. \tag{7.42}$$

Equations (7.41) and (7.42) give elemental expressions for $m_{ij}^{(2)}$ for all $i, j$.
Further $m_{ij}m_{jj}^{(2)} / m_{jj} = m_{ij}[1 + 2m_{jj}(g_{jj} - g_{bj}^{(2)})]$ leading to equations (7.39).

Finally, the variances of the first passage times $T_{ij}$, given, by equations (7.40), can be derived as var[$T_{ij}$] = var $[T_{ij} \mid X_0 = i] = m_{ij}^{(2)} - (m_{ij})^2$.

The advantage of the results of Theorems 7.5 and 7.9 is that following one matrix inversion to obtain $G_{eb}$, one can find the stationary probabilities (actually only the *b-th* row of $G_{eb}$ in this case) and the mean first passage times. In addition, computing $(G_{eb})^2$ leads to expressions for the second moments, and hence the variance of the first passage times (only the elements of $G_{eb}$ and the $b$-th rows $(G_{eb})^2$ for the first return times when $i = j$).

The efficiency of such a procedure is clear. The inaccuracies alluded earlier are reduced to a minimum with the requirement that only an accurate package to compute a single matrix inverse (of a matrix whose elements do not need to be computed in advance) is required.

## 8. Occupation Time Random Variables

Another application arises in examining the asymptotic behaviour of the number of times particular states are entered.

Given a Markov chain $\{X_n\}$, let $M_{ij}^{(n)}$ = Number of $k$ ($0 \le k \le n$) such that $X_k = j$ given $X_0 = i$. Then

$$\left[EM_{ij}^{(n)}\right] = \left[\sum_{k=0}^{n} p_{ij}^{(k)}\right] = \sum_{k=0}^{n} P^k .$$

Let $A_n = \sum_{k=0}^{n-1} P^k$ then $(I - P)A_n = I - P^n$ or $A_n(I - P) = I - P^n$, subject to $A_n\Pi = \Pi A_n = n\Pi$.

For $G$, any g-inverse of $I - P$, the standard generalised inverse method for solving such equations (Proof omitted - Exercise) leads to the following results ([Hunter, 1983):

$$A_n = \sum_{k=0}^{n-1} P^k = \begin{cases} n\Pi + (I - P)G(I - P^n), \\ n\Pi + (I - P^n)G(I - P) \end{cases}$$

$$\Rightarrow \sum_{k=0}^{n-1} P^k = n\Pi + (I - \Pi)G(I - \Pi) + o(1)E.$$

Thus $\left[EM_{ij}^{(n)}\right] = (n + 1)\Pi + (I - \Pi)G(I - \Pi) + o(1)E,$

$= (n + 1)\Pi + A^{\#} + o(1)E,$

A simple observation from this last result is that $(I - \Pi)G(I - \Pi)$ is invariant for all $G$. In fact this expression is equivalent to $A^{\#}$, the group inverse of $I - P$, giving a simple procedure for finding the group inverse. (Meyer, 1975), (Hunter, 1983).